\documentclass{birkmult}

 \newtheorem{thm}{Theorem}[section]
 
 \newtheorem{lem}[thm]{Lemma}
 
 \theoremstyle{definition}
 
 \theoremstyle{remark}
 \newtheorem{rem}[thm]{Remark}
 \newtheorem*{ex}{Example}
 \numberwithin{equation}{section}

\newcommand{\calA}{\mathcal{A}}

\newcommand{\calC}{\mathcal{C}}
\newcommand{\calD}{\mathcal{D}}

\newcommand{\calH}{\mathcal{H}}
\newcommand{\calK}{\mathcal{K}}
\newcommand{\calL}{\mathcal{L}}
\newcommand{\calM}{\mathcal{M}}

\newcommand{\calO}{\mathcal{O}}
\newcommand{\calP}{\mathcal{P}}

\newcommand{\calR}{\mathcal{R}}

\newcommand{\calV}{\mathcal{V}}
\newcommand{\calX}{\mathcal{X}}

\newcommand{\calZ}{\mathcal{Z}}

\newcommand{\calI}{\mathcal{I}}

\newcommand{\bbB}{\mathbb{B}}
\newcommand{\bbF}{\mathbb{F}}
\newcommand{\bbC}{\mathbb{C}}

\newcommand{\bbP}{\mathbb{P}}
\newcommand{\bbQ}{\mathbb{Q}}
\newcommand{\bbR}{\mathbb{R}}

\newcommand{\bbZ}{\mathbb{Z}}

\newcommand{\bfe}{\mathbf{e}}
\newcommand{\bff}{\mathbf{f}}

\newcommand{\bfk}{\mathbf{k}}

\newcommand{\la}{\langle}
\newcommand{\ra}{\rangle}

\newcommand{\s}{\textup{s}}
\newcommand{\ses}{\textup{ss}}
\newcommand{\Fl}{\text{Fl}}
\def\GL{\text{GL}}
\def\SL{{\text{SL}}}
\def\Aut{{\text{Aut}}}

\def\Pic{{\text{Pic}}}
\def\Ker{{\text{Ker}}}
\def\O{{\text{O}}}
\def\SO{\text{SO}}
\def\SU{\text{SU}}

\def\rank{{\text{rank}}}
\def\det{{\text{det}}}
\def\mod{{\text{mod}}}
\def\Tors{\text{Tors}}
\def\Mat{\text{Mat}}
\def\Sp{\text{Sp}}
\def\Hom{\text{Hom}}
\def\Alb{\text{Alb}}
\def\alb{\text{alb}}
\def\per{\text{per}}
\def\Jac{\text{Jac}}
\def\Pic{\text{Pic}}
\def\wt{\text{Wt}}
\def\bldmu{\boldsymbol{\mu}}

\begin{document}
%
%
%
%
%
%
%
%
%
\title[Moduli of K3 surfaces]
 {Moduli spaces of K3 Surfaces and Complex Ball Quotients}
\author[Igor V. Dolgachev]{Igor V. Dolgachev}

\address{%
Department of Mathematics, University of Michigan, Ann
Arnor, MI 48109,USA
}

\email{idolga@umich.edu}
\thanks{Research of the first author is partially supported by NSF grant 0245203\\
Research of the second author is partially supported by
Grant-in-Aid for Scientific Research A-14204001, Japan.}
\author{Shigeyuki Kond${\rm \bar o}$}
\address{Graduate School of Mathematics,
Nagoya University,
Nagoya 464-8602,
Japan}
\email{kondo@math.nagoya-u.ac.jp}


\date{September 1, 2005}
\dedicatory{}

\begin{abstract}
These notes are based on a series of talks given by the authors at the CIMPA Summer School on Algebraic Geometry and Hypergeometric Functions held in Istanbul in Summer of 2005.  They provide an introduction to  
a recent work  on the complex ball uniformization of the moduli spaces of Del Pezzo surfaces, K3 surfaces  and algebraic curves of lower genus. We discuss the relationship of these constructions with the Deligne-Mostow theory of periods of hypergeometric differentianl forms.  For convenience to a non-expert reader we include an introduction to the theory of periods of integrals on algebraic varieties with emphasis on abelian varieties and K3 surfaces.
\end{abstract}

\maketitle
\section{Introduction}

These notes are based on a series of talks at the CIMPA Summer School on Algebraic Geometry and Hypergeometric Functions held in Istanbul in Summer of 2005.  The topic of the talks was an introduction  to a recent work of various people on the complex ball uniformization of the moduli spaces of Del Pezzo surfaces and algebraic curves of lower genus (\cite{ACT}, \cite{DGK}, \cite{K1}-\cite{K4}, \cite{HL}). Keeping in mind the diverse background of the audimore emphasis on abelian varieties and K3 surfaces. So, an expert may start reading the notes from section 7. 

It has been  known for more than a century that a complex structure on a Riemann surface of genus $g$ is determined up to isomorphism by the period matrix $\Pi = (\int_{\gamma_j}\omega_i)$, where $(\gamma_1,\ldots,\gamma_{2g})$ is a basis of 1-homology and $(\omega_1,\ldots,\omega_g)$ is a basis of holomorphic 1-forms. It is possible to choose the bases in a such a way that the matrix $\Pi$ has the form $(Z\ |\  I_g)$, where $Z$ is a symmetric complex matrix of size $g$ with positive definite imaginary part. All such matrices are parametrized by a complex domain $\calZ_g$ in $\bbC^{g(g+1)/2}$ which is homogeneous with respect to the group $\Sp(2g,\bbR)$. In fact, it represents an example of a Hermitian symmetric space of non-compact type, a Siegel half-plane of degree $g$. A different choice of bases with the above property of the period matrix corresponds to a natural action of the group $\Gamma_g = \Sp(2g,\bbZ)$ on $\calZ_g$. In this way the moduli space $\calM_g$ of complex structures on Riemann surfaces of genus g admits a holomorphic map to the orbit space $\Gamma_g\backslash \calZ_g$ which is called the period map. The fundamental fact is the Torelli Theorem which asserts that this map is an isomorphism onto its image. This gives a moduli theoretical interpretation of some points of the orbit space. All points can be interpreted as the period matrices of principally polarized abelian varieties, i.e. $g$-dimensional complex tori equipped with an ample  line bundle whose space of holomorphic sections is one-dimensional. In this way the orbit space $\Gamma_g\backslash \calZ_g$ becomes isomorphic to the moduli space $\calA_g$ of such complex tori. The Siegel half-plane $\calZ_g$ is a Hermitian symmetric space of non-compact type III in Cartan's classification. The development of the general theory of periods of integrals on algebraic varieties in the sixties due to P.  Griffiths raised a natural question on moduli-theoretical interpretation of other Hermitian summetric spaces and their arithmetic quotients and the analogs of the Torelli Theorem.  For the spaces of classical types  I-IV, this can be achieved by embedding any type space into a Siegel half-plane, and introducing the moduli space of abelian varieties with some additional structure (like a level, complex multiplication or some tensor form on cohomology). All of this becomes a part of the fancy theory of Shimura varieties. However, a more explicit interpretation remains to be searched for. For the type IV Hermitian symmetric space of dimension $\le 19$ such an interpretation had been found in terms of moduli of complex algebraic surfaces of type K3. The fundamental result of I. Piatetsky-Shapiro and I. Shafarevich in the seventies gives an analog of the Torelli Theorem for polarized algebraic K3 surfaces. Although this gives only a realization of the 19-dimensional type IV space by a very special arithmetic group depending only on the degree of the polarization, a modified version of the polarization structure due to V. Nikulin, allows one to extend this construction to type IV space of any dimension $d\le 19$ with a variety of arithmetric groups realized as certain orthogonal groups of integral quadratic forms of signature $(2,d)$. Recently, some arithmetic quotients of type IV domain of dimension 20 have been realized as periods of holomorphic symplectic manifolds of dimension 4.

A complex ball is an example of Hermitian symmetric space of type I. Some of its arithmetic quotients have a realization as periods of hypergeometric differential forms  via the Deligne-Mostow theory. We refer for the details to Looijenga's article in the same volume \cite{Lo}. A hypergeometric differential form has an interpretation as a holomorphic 1-form on a certain algebraic curve on which a cyclic group acts by automorphisms and the form is transformed according to a character of this group. In section 7 of our paper we discuss in more general setting  the theory of what we call  eigenperiods of algebraic varieties. 

The periods of hypergeometric functions allows one to realize some complex ball quotients as the moduli space of weighted semi-stable ordered point sets in projective line modulo projective equivalence. In some cases these moduli spaces are isomorphic to moduli spaces of other structures. For example, via the period map of algebraic curves, the hypergeometric  complex ball quotients are mapped onto a subvariety  of  $\calA_g$ parametrizing principally polarized abelian varieties with a certain cyclic group action. For example, the moduli space of 6 points with equal weights  defines a curve of genus 4 admitting a cyclic triple cover of the projective line ramified at the 6 points. Its Jacobian variety is an abelian variety of dimension 4 with a cyclic group of order 3 acting by automorphism of certain type. The moduli space of such abelian varieties is a complex ball quotient.  Another example is the moduli space of equally marked sets of 5 points which leads to the moduli space of marked Del Pezzo surfaces of  degree 4 \cite{K3}. Other examples relating the arithmetic complex ball quotients arising in the Deligne-Mostow theory to moduli spaces of Del Pezzo surfaces were found in \cite{MT}, \cite {HL}. It turns out that all of these examples are intimately related to the moduli space of K3 surfaces with specail structure of its Picard group of algebraic cycles and an action of a cyclic group. In section 11 we develop a general theory of such moduli spaces. In section 12 we briefly discuss all the known examples of moduli spaces of Del Pezzo surfaces and curves of low genus which are isomorphic to the moduli space of such structres on K3 surfaces, and via this isomorphism admit a complex ball uniformization by an arithmetic group. Some of these examples arise from the Deligne-Mostow theory.  In the last section 12 we discuss a relationship between  Deligne-Mostow arithmetic complex ball quotients and the  moduli spaces of K3 surfaces with a cyclic group action. 

We would like to thank the organizers of the Summer School, and especially Professor A. Muhammed Uludag for the hospitality and for providing a stimulating and pleasant audience for our lectures. 

\section{Introduction to Hodge theory}

In this and next four sections we give a brief introduction to the Hodge theory of periods of integrals on algebraic varieties. We refer for details to  \cite{CMP}, \cite{GH}, \cite{Vo}. 

Let $M$ be a smooth compact oriented connected manifold of even dimension $2n$. Its 
cohomology $H^*(M,\bbZ)$ is a graded algebra over $\bbZ$ 
with multiplication defined by the cup-product
$$\cup:H^k(M,\bbZ)\times H^l(M,\bbZ) \to H^{k+l}(M,\bbZ)$$
satisfying $x\cup y = (-1)^{kl} y\cup x$. In particular, the restriction 
of $\cup$ to the middle-dimensional cohomology $H^n(M,\bbZ)$ is a $\bbZ$-bilinear 
form
$$b_M:H^n(M,\bbZ)\times H^n(M,\bbZ) \to H^{2n}(M,\bbZ) \cong \bbZ,$$
where the latter isomorphism is defined by using the fundamental class $[M]$ of $M$. The 
Poincar\'e duality asserts that this bilinear form is a perfect pairing, modulo torsion. It is also symmetric if $n$ is even, and skew-symmetric if $n$ is odd.

Recall that the cohomology $H^*(M,\bbR) = H^*(M,\bbZ)\otimes \bbR$ can be
 computed by using the De Rham theorem:
$$H^*(M,\bbR) \cong H^*(0\to A^0(M)\overset{d}{\to} A^1(M)\overset{d}{\to} A^2(M)\to \ldots),$$
where  $A^k(M)$ is the space of smooth differential $k$-forms on $M$. The cup-product is defined by 
$$[\alpha]\cup [\beta](\gamma)= \int_\gamma\alpha\wedge\beta, \ \gamma\in
 H_*(M,\bbR),$$
where we consider cohomology as linear functions on homology. In particular, the bilinear form $b_M$ is an inner product on $H^n(M,\bbR)$ defined by the formula
\begin{equation}\label{cup}
b_M([\alpha],[\beta]) = \int_M\alpha\wedge \beta.\end{equation}
The same is true for the cohomology $H^*(M,\bbC) = H^*(M,\bbR)\otimes \bbC$ if we replace $A^k(M)$ with complex valued smooth $k$-forms. Now suppose $M$ is the underlying differentiable structure of a compex manifold $X$. Then local coordinates $t_1,\ldots,t_{2n}$ can be expressed in terms of complex coordinates
 $z_1,\ldots,z_n$ and its conjugates $\bar{z}_1,\ldots,\bar{z}_n$. This allows us to express locally a smooth $k$-form as a sum of   forms of type $(p,q)$:
\begin{equation}\label{form}\omega = \sum a_{i_1\ldots i_p;j_1\ldots j_q}
(z,\bar{z})dz_{i_1}\wedge\ldots\wedge dz_{i_p}\wedge d\bar{z}_{j_1}\wedge\ldots \wedge d\bar{z}_{q}, \quad p+q = k,
\end{equation}
where $a_{i_1\ldots i_p;j_1\ldots j_q}(z,\bar{z})$ are smooth complex-valued functions. This gives a decomposition 
$$A^k(X) = \bigoplus_{p+q=k} A^{pq}(X)$$
and $d= d'+d''$, where $d'$ (resp. $d''$) is the derivation operator with respect to the variables $z_i$ (resp.  $\bar{z}_i)$). The Dolbeaut Theorem gives an isomorphism
$$H^q(X,\Omega_X^p) \cong  H^q(0\overset{d''}{\to} A^{p,0}\overset{d''}{\to} A^{p,1}\overset{d''}{\to} \ldots),$$
 where $\Omega_X^p$ is the sheaf of holomorphic $p$-forms, i.e. forms from
 \eqref{form} of type $(p,0)$, where the coefficients are holomorphic functions. The Dolbeaut Theorem generalizes a well-known fact that a smooth function $f(z,\bar{z})$ of a complex variable is holomorphic if and only if it satisfies the equation $\frac{\partial f(z)}{\partial \bar{z}} = 0$.

A structure of a complex manifold $X$ on  a smooth manifold gives a decomposition of the complexified tangent bundle $T_M\otimes\bbC$ into a holomorphic and anti-holomorphic part with local basis $\bigl(\frac{\partial}{\partial z_i}\bigr)$ and $\bigl(\frac{\partial}{\partial \bar{z}_i}\bigr)$, respectively. We denote the holomorphic part by $T_X$. An additional structure of a hermitian complex manifold on $X$ is given by a choice of a holomorphically varying structures of  a hermitian vector space on tangent spaces $T_{X,x}$ defined by a tensor
\begin{equation}\label{metric}
ds^2 = \sum h_{ij}(z)dz_i\otimes d\bar{z}_j.
\end{equation} 
It allows one to define the adjoint operator $\delta''$ of $d''$, the Laplace operator $\Delta'' = \delta''d''+d''\delta''$, and the notion of a {\it harmonic form} of type $(p,q)$ (an element of the kernel of the Laplace operator). One shows that each $d''$-closed form of type $(p,q)$ is $d''$-cohomologous to a unique harmonic form of type $(p,q)$. In particular, there is a canonical isomorphism of vector spaces
 $$\calH^{pq}\otimes \bbC \cong H^q(X,\Omega_X^p),$$
where $\calH^{pq}$ is the space of harmonic forms of type $(p,q)$. 
 On the other hand, a hermitian form defines a structure of a Riemannian manifold on $M$. The latter defines the adjoint operator $\delta$ of $d$, the Lapalace operator $\Delta = \delta d+d\delta$, and the space of harmonic forms $\calH^k$. Each cohomology class has a unique representative by a harmonic form, i.e. there is a canonical isomorphism of vector  spaces
 $$H^k(M,\bbR) = \calH^k.$$
 A fundamental fact proved by Hodge asserts that if a hermitian metric on $X$ is a {\it K\"ahler metric}, i.e. the $(1,1)$-form  
 \begin{equation}\label{kahler}
 \omega =  \frac{i}{2}\sum h_{ij}(z)dz_i\wedge d\bar{z}_j
 \end{equation}
 is closed (in this case the 2-form is  called the {\it K\"ahler form}),  the Laplace operators
 $\Delta''$ and $\Delta$ extended to  $H^*(X,\bbC)$ coincide up to a constant factor. This shows that 
 $$\sum_{p+q=k} \calH^{pq} = \calH^k\otimes \bbC,$$ and we have a {\it Hodge decomposition}
\begin{equation}\label{hdecomp}
H^k(X,\bbC) = \bigoplus_{p+q=k} H^{pq}(X),
\end{equation}
where 
$$H^{pq}(X)  \cong H^q(X,\Omega_X^p).$$
The Hodge decomposition satisfies the following properties
\begin{itemize}
\item[(HD1)] The decomposition \eqref{hdecomp} does not depend on the choice of a K\"ahler metric;
\item[(HD2)] $H^{pq} = \overline{H^{qp}}$, where the bar denotes the complex conjugation. In particular $h^{pq}(X) = h^{qp}(X)$, where
$$h^{pq}(X) = \dim_\bbC H^{pq}(X);$$
\item[(HD3)]  The bilinear form $Q:H^k(X,\bbR)\times H^k(X,\bbR)\to \bbR$ defined by 
the cup-product $(\phi,\psi)\mapsto \phi\cup \psi\cup [\omega]^{n-k}$ is symmetric when $k$ is even and skew-symmetric otherwise. It satisfies
$$Q(x,y) = 0, \quad x\in H^{pq}, y\in H^{p'q'}, p\ne q'.$$
\item[(HD4)] Let $H_{prim}^k(X,\bbR)$ denote the orthogonal complement of $[\omega]^{n-\frac{k}{2}}$ with respect to $Q$ (for $k$ odd $H^k(X,\bbR)_{prim}: = H^k(X,\bbR))$. For any nonzero $x\in H^{pq}_{prim}:=H^{pq}\cap H_{prim}^k(X,\bbR)$,
$$i^{p-q}(-1)^{k(k-1)/2}Q(x,\bar{x}) > 0.$$
 \end{itemize}

Using property (HD3)  one can compute the signature $I(M)$  of the cup-product on $H^n(M,\bbR)$ if $n$ is even and $M$ admits a K\"ahler structure:
\begin{equation}
I(M) =  \sum_{p\equiv q\mod 2}(-1)^ph_1^{pq},
\end{equation}
where $h_1^{pq} = h^{pq}$ if $p\ne q$ and $h_1^{pq} = h^{pq}-1$ otherwise. In particular, the Sylvester signature $(t_+,t_-)$ (recall that, by definition, $I(M) = t_+-t_-$) of the quadratic form on $H^n(M,\bbR)$ when $n$ is even, is given by 
\begin{equation}\label{index}t_+ = \frac{1}{2}(b_n(M)+I(M)+1), \quad t_- = \frac{1}{2}(b_n(M)-I(M)
+1),\end{equation}
where $b_i(M) = \dim_\bbR H^i(M,\bbR)$ are the Betti numbers of $M$.

Define the subspace $F^p$ of $H^k(X, \bbC)$ by
$$F^p = \sum_{p'\ge p} H^{p'q}(X),$$
so that $F^0 = H^k(X,\bbC), F^k  = H^{k,0}, F^{p} = \{0\}, p > k$. This defines a flag $(F^p)$ of linear subspaces 
$$0 \subset F^k \subset F^{k-1} \subset \ldots \subset F^0 = 
H^k(X,\bbC).$$
Note that the Hodge decomposition can be reconstructed from this flag by using property (HD2)
\begin{equation}\label{rec}
H^{pq}(X) = \{x\in F^{p}: Q(x,\bar{y}) = 0, \ \forall y\in F^{p+1}\}.
\end{equation}
Assume that the K\"ahler form \eqref{kahler} is a {\it Hodge form}, i.e. its cohomology class $[\omega]$  belongs to $H^2(X,\bbZ)$.  By a theorem of Kodaira this implies that $X$ is isomorphic to a complex projective algebraic variety. It admits a projective embedding such that the cohomology class of a hyperplane section is equal to some positive multiple of $[\omega]$. The Hodge decomposition in this case has an additional property that the form $Q$ takes integer values on the image of $H^k(X,\bbZ)$ in $H^k(X,\bbR)$.

The flag $(F^p)$ defined by a Hodge decomposition is an invariant of a complex structure on $M$, and a bold conjecture (not true in general) is that it completely determines the complex structure on $M$ up to isomorphism. In fact, this was shown in the 19th century for Riemann surfaces, i.e. complex manifolds of dimension 1.

One can define an {\it abstract Hodge structure of weight $k$} (AHS) on a real vector space $V$  to be a decomposition into direct sum of complex subspaces
$$V_\bbC = \bigoplus_{p+q= k} V^{pq}$$
such that $\overline{V^{pq}} = V^{qp}$. A {\it polarization} of AHS on $V$ is a a non-degenerate bilinear form $Q$ on $V$ which is symmetric if $k$ is even, and skew-symmetric otherwise. It satisfies  the condition (HD3) and (HD4) from above, where $H^{pq}_{prim}$ is replaced with $V^{pq}$.
 
 An {\it integral structure} of a AHS is a free abelian subgroup $\Lambda\subset V$ of rank equal to $\dim V$ (a {\it lattice}) such that 
$Q(\Lambda\times \Lambda)\subset \bbZ$. One can always find an integral structure by taking $\Lambda$ to be the $\bbZ$-span of a {\it standard basis} of $Q$ for which the matrix of $Q$ is equal to a matrix
\begin{equation}\label{ipq}
I(a,b) = \begin{pmatrix}I_a&0_{a,b}\\
0_{b,a}&-I_{b}\end{pmatrix}.\end{equation}
if $k$ is even, and the matrix
\begin{equation}\label{ssb}
J = \begin{pmatrix}0_{m}&I_m\\
-I_m&0_{m}\end{pmatrix}.\end{equation}
if $k$ is odd.

 A Hodge structure on cohomology $V = H_{prim}^k(X,\bbR)$  is an example of an AHS of weight $k$. If $[\omega]$ is a Hodge class, then the cohomology admits an integral structure with respect to the lattice $\Lambda$ equal to the intersection of the image of $H^k(X,\bbZ)$ in $H^k(X,\bbR)$ with $H_{prim}^k(X,\bbR)$. 

\section{The period map}\label{}

Let $(F^p)$ be the flag of subspaces of $V_\bbC$ defined by a polarized  AHS of weight $k$ on a vector space $V$. Let
$$f_p = \sum_{p'\ge p}h^{p'q}(X) =\dim F^p, \quad \bff = (f_0,\ldots,f_k).$$
 Let $\Fl(\bff,V_\bbC)$ be the  variety  of flags of linear subspaces $F^p$ of dimensions $f_p$, $p=0,...,k$.
It is a closed algebraic subvariety of the product of the Grassmann varieties  $G(f_p,V_\bbC)$. 
 A polarized  AHS of weight $k$ defines a point $(F^p)$ in $\Fl(\bff,V_\bbC)$. It satisfies the following conditions 
\begin{itemize}
\item[(i)] $V_\bbC = F^p\oplus \overline{F^{k-p+1}}$;
\item[(ii)] $Q(F^p,F^{k-p+1}) = 0$;
\item[(iii)] $(-1)^{k(k-1)/2}Q(Cx,\bar{x}) > 0$, where $C$ acts on $H^{pq}$ as multiplication by  
$i^{p-q}$.
\end{itemize}
The subset of flags in  $\Fl(\bff,V_\bbC)$ satisfying the previous conditions is denoted by 
$\calD_\bff(V,Q)$ and is called the {\it period space} of $(V,Q)$ of type $\bff$.

Fix a standard basis in $V$ with respect to $Q$ to identify $V$ with the space $\bbR^{r}$, where $r = f_0$, and denote by
$A$ the corresponding matrix of $Q$.
Let  $F^p$ be the column space of a complex matrix $\Pi_p$ of size $r\times f_p$.  We assume that the first $f_{p+1}$ columns of $\Pi_{p}$  form the matrix $\Pi_{p+1}$. Then any flag in 
$\calD_\bff:=\calD_\bff(\bbR^{r},A)$ is described by a set of matrices $\Pi_p$ satisfying  the following conditions:
\begin{itemize}
\item[(PM1)] $\det (\Pi_p|\overline{\Pi_{k-p+1}}) \ne 0$;
\item[(PM2)] ${}^t\Pi_p\cdot A\cdot \Pi_{k-p+1} = 0$;
\item[(PM3)] $(-1)^{k(k-1)/2}{}^t\Pi_p\cdot A\cdot \overline{\Pi_{k-p+1}}\cdot K> 0$, 
where $K$ is a diagonal matrix with $\pm i$ at the diagonal representing the operator $C$. 
\end{itemize}
Note that condition (ii) is a closed algebraic  condition and defines a closed algebraic subvariety of $\Fl(\bff,V_\bbC)$. Other conditions are  open conditions in complex topology. Thus the collections of matrices $(\Pi_p)$   has a natural structure of a complex manifold.  Two  collections of matrices $(\Pi_p)$ and $(\Pi_p')$ define the same point in the period space if and only if there exists an invertible complex matrix $X$ of size $r$  such that 
$\Pi_0' = \Pi_0\cdot X$. The matrix $X$  belongs to a subgroup $P_\bff$  of 
$\GL(r,\bbC)$ preserving the flag of the subspaces $F^p$ generated by the first $f_p$ unit vectors 
$\bfe_i$. The flag variety $\Fl(\bff,V_\bbC)$ is isomorphic to the homogeneous space 
$$\Fl(\bff,V_\bbC) \cong  \GL(r,\bbC)/P_\bff.$$
The period space   $\calD_{\bff}$ is an open subset (in complex topology) of a  closed algebraic subvariety  $\check{\calD}_{\bff}$ of  $\Fl(\bff,V_\bbC)$ defined by condition (ii). It is known  that the group $G_\bbC = \Aut(\bbC^r,Q_0)$ acts transitively on $\check{\calD}_{\bff}$ with isotropy subgroup $P = G_\bbC\cap P_\bff$ so that 
$$\check{\calD}_{\bff}\cong G_\bbC/P$$
is a projective homogeneous variety. The group $G_\bbR= \Aut(\bbR^r,Q_0)$ acts transitively on 
$\calD_{\bff}$ with a compact isotropy subgroup $K$ so that
$$\calD_{\bff} \cong G_\bbR/K$$
is a complex non-compact homogenous space. 

In the case when the AHS is the  Hodge structure on cohomology $H^k(X,\bbR)$ of a K\"ahler manifold $X$, the matrices $\Pi_p$ are called the {\it period matrices}. If $(\gamma_1,\ldots,\gamma_{f_0})$ is a basis in $H_k(X,\bbZ)/\Tors$ such that the dual basis $(\gamma_1^*,\ldots,\gamma_{f_0}^*)$ is a standard basis of the polarization form defined by a choice of a K\"ahler structure, then 
$$\Pi_p = (\int_{\gamma_i}\omega_j),$$
where $(\omega_1,\ldots,\omega_{f_p})$ is a basis of $F^p$ represented by differential $k$-forms.

Now, suppose  we have a {\it family} of compact connected complex manifolds. It is a holomorphic smooth map 
$f:\calX \to T$ of complex manifolds  with connected base $T$ and with its fibre $X_t = f^{-1}(t)$. For any point $t\in T$ we have the real vector space $V_t = H^k(X_t,\bbR)$ equipped with a Hodge structure. We also fix a K\"ahler class  $[\omega]$ on $\calX$ whose restriction to each $X_t$ defines a polarization $Q_t$ of the Hodge structure on $V_t$.  One can prove (and this is not trivial!) that the Hodge numbers $h^{pq}(X_t)$ do not depend on $t$. Fix  an isomorphism 
\begin{equation}\label{mark}
\phi_t:(V,Q_0) \to (H^k(X_t,\bbR),\cup)
\end{equation}
called a {\it $k$th marking} of $X_t$. Then pre-image of the Hodge flag $(F_t^p)$ is a  $Q_0$-polarized  AHS of weight $k$ on $(V,Q_0)$. Let $\calD_\bff$ be the period space of $(V,Q_0)$ of type $\bff =(f_p)$, where $f_p = \dim F_t^p$. We have a ``multi-valued '' map 
$$\phi:T\to \calD, \quad t\mapsto \bigl(\phi^{-1}(F_t^p)\bigr).$$
called the {\it multi-valued period map} associated to $f$. According to a theorem of Griffiths this map is a multi-valued holomorphic map. 

To make the period map one-valued one we need to fix a standard basis in $(V_t,Q_t)$ for each $t$ which depends  holomorphically on $t$. This is not possible in general. The vector spaces   $(V_t)_{t\in T}$ form a real {\it local coefficient system} $\calV$ on $T$, i.e. a real vector bundle of rank $f_0$ over $T$ whose transition functions  are matrices with constant entries.  Any local coefficient system is determined by its {\it monodromy representation}. If we fix a point $t_0\in T$ and let $V = V_{t_0}$, then for any continuous loop $\gamma:[0,1]\to T$ with $\gamma(0) = \gamma(1) = t_0$, the pull-back 
$\gamma^*(\calV)$ of $\calV$ to $[0,1]$ is a trivial  local coefficient system which defines a linear self-map 
$a_\gamma:V = \gamma^*(\calV)_0 \to V = \gamma^*(\calV)_1$. It depends only on the homotopy class $[\gamma]$ of $\gamma$. The map $[\gamma]\mapsto a_\gamma$ defines a  homomorphism of groups
$$a:\pi_1(T,t_0) \to \GL(V), $$
called the {\it monodromy representation}. The pull-back of $\calV$ to the universal covering 
$\tilde{T}$ is isomorphic to the trivial local coefficient system $\tilde{T}\times V$ and 
$$\calV \cong \tilde {T}\times V/\pi_1(T,t_0),$$
where $\pi_1(T,t_0)$ acts by the formula $[\gamma]:(z,v) =([\gamma]\cdot z, a([\gamma])(v))$. Here the action of $\pi_1(T,t_0)$ on $\tilde{T}$ is the usual action by deck transformations.

One can show that the monodromy representation preserves the polarization form $Q_{t_0}$ on $V$ so that the image $\Gamma(f)$ of the monodromy representation lies in $\Aut(V,Q)$. This image is called the {\it monodromy group} of $f$. Of course, if $T$ happens to be simply-connected, say by restriction to a small neighborhood of a point $t_0$, the monodromy representation is trivial, and we can define a one-valued period map. 

Let  $\pi:\tilde{T}\to T$ be the universal covering map. The second projection of the fibred product 
$\tilde{\calX} = \calX\times_T\tilde{T}$ of complex manifolds defines a holomorphic map
$\tilde{f}:\tilde{\calX} \to \tilde{T}.$
Its fibre over a point $\tilde{t}$ is isomorphic to the fibre $X_t$ over the point $t = \pi(\tilde{t})$. Now we can define the one-valued holomorphic map
\begin{equation}\label{period}
\tilde{\phi}:\tilde{T} \to \calD_\bff.
\end{equation}

Assume that the monodromy group $\Gamma(f)$ is a discrete subgroup of $\Aut(V,Q)$ with respect to its topology of a real Lie group. Since $\Aut(V,Q)$ acts on $\calD_\bff$ with compact isotropy subgroups, the isotropy subgroups of $\Gamma(f)$ are finite. This implies that the orbit space $\Gamma(f)\backslash \calD_\bff$ is a complex variety with only quotient singularities
and the natural projection to the orbit space is a holomorphic map. The period map $\tilde{\phi}$ descends to a one-valued holomorphic map
\begin{equation}\label{periodq}
\bar{\phi}: T \to \Gamma(f)\backslash \calD_\bff.
\end{equation}

The discreteness condition for the monodromy group is always satisfied if $f:\calX\to T$ is a family of projective algebraic varieties and the K\"ahler form on each $X_t$ defining the polarization $Q_t$ is a Hodge form. In this case the monodromy group is a subgroup of $G_\Lambda = \Aut(\Lambda,Q|\Lambda)$, where 
$\Lambda$ is the image of $H^k(X_{t_0},\bbZ)$ in $V$.  Thus for any such family there is a holomorphic map
\begin{equation}
\bar{\phi}: T \to G_\Lambda\backslash \calD_\bff.
\end{equation}

One says that a K\"ahler manifold $X$ satisfies an  {\it Infinitesimal Torelli Theorem} if for any family of varieties $f:\calX \to T$ as above with $X\cong X_{t_0}$ for some $t_0\in T$ the period map \eqref{period} is  an isomorphic embedding of some analytic neighborhood of $t_0$ (or, equivalently, the differential of the map at $t_0$ is injective).  We say that $X$ satisfies a {\it Global Torelli Theorem} if  for any two points  $t,t'\in T$ with the same image in $\Gamma(f)\backslash\calD_\bff$ there is an isomorphism of complex manifolds $\phi:X_t\to X_{t'}$ such that $f^*([F^p_{t'}]) = [F^p_t]$.  

Fix a smooth manifold $M$  underlying some complex  manifold. Consider the set of isomorphism classes of complex structures on $M$. A {\it moduli problem} is the problem of putting on this set a structure of an analytic space (or a complex variety) $\calM$  such that  
 for any holomorphic map $f:Y\to T$ of analytic spaces whose fibres are complex manifolds diffeomorphic  to $M$, the map $T\to \calM$ which assigns to $t\in T$ the isomorphism class of $f^{-1}(t)$ defines a holomorphic map $\phi:T\to \calM$. Or, equivalently, for any $f:Y\to T$ there exists a holomorphic map $\phi:T\to \calM$ such that  the fibres which are mapped to the same point must be isomorphic complex manifolds.  We require also that for any $\calM'$ with the same property  there exists a holomorphic map $s:\calM \to \calM'$ such that $\phi' = s\circ \phi$. If $\calM$ exists it is called a {\it coarse moduli space} of complex structures on $M$.

 If additionally, there exists a holomorphic map
$u:\calX \to \calM$ such that any $f$ as above is obtained via a unique map $\phi:T\to \calM$ by taking the fibred product $\calX\times_\calM T\to T$. In this case $\calM$ is called a {\it fine moduli space}, and $u$ the {\it universal family}. Note that, for any point $m\in \calM$, the isomorphism class of the fibre $u^{-1}(m)$ is equal to $m$. Fine moduli spaces rarely exist unless we put some additional data, for example a marking on cohomology as in \eqref{mark}. 

Similar definitions one can give for the moduli space of structures of an algebraic variety, or for polarized K\"ahler manifolds or polarized projective algebraic varieties. For the latter, one fixes a cohomology class in $H^2(M,\bbR)$  and consider the set of complex structures on $M$ such that this cohomology class can be represented by a  K\"ahler (or a Hodge) form. 

Suppose a coarse moduli of polarized algebraic manifolds of diffeomorphism type $M$ exists and the Global Torelli Theorem holds for the corresponding complex manifolds. Then, by definition of a coarse moduli space there must be a holomorphic map
$$\text{per}:\calM \to  G_\Lambda\backslash \calD_\bff.$$
A {\it Local  Torelli Theorem} is the assertion that this map is a local isomorphism which, together with Global Torelli Theorem, will assert that the map is an embedding of complex varieties.  Note that the Local Torelli Theorem implies the Infinitesimal Torelli Thorem but the converse is not true in general.

\section{Hodge structures of weight 1}\label{}

An abstract Hodge structure of weight 1 on a real vector space $V$
is a decomposition into direct sum of  complex linear subspaces
\begin{equation}\label{hodec}
V_\bbC = V^{10}\oplus V^{01}.
\end{equation}
such that 
\begin{itemize}
\item [(i)] $\overline{V^{10}} = V^{01}$.
\end{itemize}

A polarization on AHS of weight 1 is a non-degenerate skew-symmetric form $Q$ on $V$ (a symplectic form)  such that 
\begin{itemize}
\item[(ii)] $Q|V^{10} = 0, Q|V^{01} = 0$;
\item[(iii)] $iQ(x,\bar{x}) > 0, \quad \forall x\in V^{10}\setminus \{0\}.$
\end{itemize}
Here, as always,  we extend $Q$ to the complexification by linearity. 

Consider a Hodge structure of weight 1 on $V$. Let $W = V^{01}$.  The composition $V\to W$ of the natural inclusion map $V\hookrightarrow V_\bbC$ and the projection to $W$ with respect to the decomposition \eqref{hodec} is a $\bbR$-linear isomorphism. This allows to transfer the structure of a complex space on $W$ to $V$. Recall that a complex structure on a real vector space $V$ is defined by a linear operator $I$ satisfying $I^2 = -1_V$ by setting
$$(a+bi)\cdot v = av+bI(v). $$
 We have 
$$V = \{w+\bar{w}, w\in W\}$$
and, for any $v=w+\bar{w}$ we  define $I(v)$ to be the unique vector in $V$ such that
$$I(v) = iw-i\bar{w}.$$
In particular, the $\bbR$-linear isomorphism $V\to W$ is defined by $v\mapsto \frac{1}{2}(v-iI(v))$. We will often identify $V$ with $W$ by means of this isomorphism. Conversely, a complex structure $I$ on $V$ defines a decomposition \eqref{hodec}, where $W$ (resp. $\overline{W}$) is the eigensubspace of $I$ extended to $V_\bbC$ by linearity with eigenvalue $i =\sqrt{-1}$ (resp. $-i$). Thus we obtain a bijective correspondence beween AHS of weight 1 on $V$  and complex structures on $V$. 
Note that replacing the Hodge structure by switching $V^{10}$ with $V^{01}$ changes the complex structure to the conjugate one defined by the operator $-I$.

Now let us see the meaning of a polarized AHS of weight 1. The polarization form $Q$ makes the pair  $(V,Q)$  a real symplectic vector space. 

 A complex structure $I$ on $V$ is called a {\it positive complex structure} with respect to $Q$ if $Q(v,I(v'))$ defines a symmetric positive definite bilinear form on $V$. It follows from the symmetry condition that the operator $I$ is an isometry of the symplectic space $(V,Q)$.

Consider a polarized Hodge structure on $V$ of weight 1 and let $I$ be the complex structure $I$ on $V$  determined by the subspace $W = V^{01}$.  Extending by linearity, we will consider $Q$ as  a skew-symmetric form on $V_\bbC$.  Let  $w = v-iI(v), w' =v'-iI(v')\in W$. We have
\begin{equation}
-iQ(w,\bar{w'}) = Q(v,I(v'))+iQ(v,v').
\end{equation} 
By definition of a polarized AHS, $Q(v,I(v)) = -iQ(w,\bar{w}) > 0$ for any $w\ne 0$. This shows that the complex structure $I$ is positive with respect to $Q$. 

Conversely, suppose $I$ is a positive complex structure  on $(V,Q)$. Let $V^{01}$ be the $i$-eigensubspace of $V_\bbC$ of the operator $I$. Since 
$Q(w,I(w')) = iQ(w,w')$ is a symmetric and also a skew-symmetric bilinear form on $W$, it must be zero. Hence $W$ (and, similarly, $V^{10} = \overline{W}$) is an isotropic subspace of $Q$.  This checks property (ii) of the Hodge structure. We have 
$0<-iQ(v-iI(v),v+iI(v)) = 2Q(v,v)$. This checks property (iii). Thus we have proved the following.

\begin{lem} There is a natural bijection between the set of Hodge structures of weight 1 on $V$ with polarization form $Q$ and the set of positive complex structures on $V$ with respect to $Q$.
\end{lem}

\begin{ex} Let $V = \bbR^{2g}$ with standard basis formed by the unit vectors $\bfe_k$. Define the complex structure $I$ by 
$$I(\bfe_k) = \bfe_{k+g}, \ I(\bfe_{g+k}) = -\bfe_k, \quad k = 1,\ldots,g.$$
The space $V_\bbC = \bbC^{2g}$ decomposes into the direct sum of the $\pm i$-eigensubspaces $V_\pm$ of $I$, where
$V_\pm$ is spanned by the vectors $\bfe_k\mp i\bfe_{k+g}.$
Let $Q$ be a skew-symmetric bilinear form defined by the condition 
$$Q(\bfe_{k},\bfe_{k+g}) = 1,\ 
Q(\bfe_{k},\bfe_{k'}) = 0 \ \text{if $|k-k'| \ne g$}$$ (the {\it standard symplectic form} on $\bbR^{2g}$). Then the matrix of $Q(v,I(v'))$ in the standard basis is the identity matrix, so $Q(v,I(v'))$ is symmetric and positive definite.
\end{ex}

Recall that a {\it hermitian form} on a complex vector space $E$ is a $\bbR$-bilinear form $H:E\times E\to \bbC$ which is $\bbC$-linear in the first argument and satisfies 
\begin{equation}\label{herm}
H(x,y) = \overline{H(y,x)}.
\end{equation}
One can view a hermitian form as a $\bbC$-bilinear map $E\times \bar{E}\to \bbC$ satisfying \eqref{herm}. Here $\bar{E}$ is the same real space as $E$ with conjugate complex structure. The restriction of $H$ to the diagonal is a real-valued quadratic form on the real space $E$, and the signature of $H$ is the signature of this quadratic form. In particular, we can speak about positive definite hermitian forms. 

\begin{lem}\label{lem1} The formula
\begin{equation}\label{herm2}
H(x,y) = -iQ(x,\bar{y})
\end{equation}
defines a hermitian form on $V_\bbC$ of signature $(g,g)$.
\end{lem}

\begin{proof} Obviously, formula \eqref{herm} defines a bilinear form on $V_\bbC\times \overline{V_\bbC}.$ Write $x,y\in V_\bbC$ in the form 
$x= v+iw, y = v'+iw'$ for some $v,w,v',w'\in V$.  Then
$$H(x,y) = -iQ(v+iw,v'-iw') = -i(Q(v,v')+Q(w,w'))-Q(w,v')+Q(v,w')$$ 
$$= \overline{H(y,x)}.$$
Let  $e_1,\ldots,e_{2g}$ be a standard symplectic basis in $V$. Let 
$f_k = e_k-ie_{k+g}, \bar{f}_k = e_k+ie_{k+g}, k = 1,\ldots, g.$ These vectors form a basis of $V_\bbC$. We have
$$H(f_k,f_l) = -iQ(f_k,\bar{f}_l) = -iQ(e_k-ie_{k+g},e_l+ie_{l+g}) = Q(e_k,e_{l+g})-Q(e_{k+g},e_l) = 2$$
if $k = l$ and $0$ otherwise.
This shows that for any nonzero $x=\sum a_kf_k$, we have
$$H(x,x) = 2\sum a_k\bar{a}_k > 0.$$
Thus the restriction of $H$ to the span $W$ of the $f_k$'s is positive definite. Similarly, we check that the restriction of $H$ to  the span of the $\bar{f}_k$'s is negative definite and two subspaces $W$ and $\overline{W}$ are orthogonal with respect to $H$. This shows that $H$ is of signature $(g,g)$.
\end{proof}

Note that the skew-symmetric form $Q$ on $V_\bbC$ is reconstructed from
$H$ by the formula
\begin{equation}\label{assoc}
Q(x,y) = iH(x,\bar{y}).\end{equation}
Let $G(g,V_\bbC)$ be the Grassmann variety of $g$-dimensional subspaces of $V_\bbC$. Let $H$ be a hermitian form on $V_\bbC$ of signature $(g,g)$ and let $Q$ be a skew-symmetric form on $V_\bbC$ associated to $H$ by the formula \eqref{assoc}.
Set
\begin{equation}\label{gras}
G(g,V_\bbC)_H = \{W\in G(g,V_\bbC): Q|W = 0, H|W > 0\}.
\end{equation}
We see that a Hodge structure of weight 1 on $V$ defines a point $W = V^{01}$ in $G(g,V_\bbC)_H$, where $H$ is defined by the formula \eqref{herm}.  Conversely, for any $W\in G(g,V_\bbC)_H$, set $V^{01} = W, V^{10} =\overline{W}$. Since $Q|W = 0$, the formula \eqref{assoc} implies that $W$ and $\overline{W}$ are orthogonal to each other, hence we have a direct sum decomposition \eqref{hodec}. We take $Q$ to be defined by \eqref{assoc}. Properties (ii) and (ii) are obviously satisfied. 

The hermitian form $H$ defined in Lemma \ref{lem1} will be called the {\it hermitian form associated} to $Q$.

We have proved the following.

\begin{thm} Let $(V,Q)$ be a real symplectic space. There is a natural bijection between Hodge structures on $V$ of weight 1 with polarization form $Q$ and points in $G(g,V_\bbC)_H$, where $H$ is the associated hermitian form of $Q$.
\end{thm}

By choosing a standard symplectic basis in $V$, $G(g,V_\bbC)_H$ can be described as a set of complex $2g\times g$-matrices satisfying conditions (PM1)-(PM3) with $p=1$:
\begin{equation}\label{riem1}
{}^t\Pi\cdot J\cdot \Pi = 0
\end{equation}
\begin{equation}\label{riem2}
-i{}^t\Pi\cdot J\cdot\bar{\Pi} >
0.
\end{equation}
Two such matrices $\Pi'$ and $\Pi$ define the same AHS if and only if there exists an invertible matrix $X$ such that $\Pi' = \Pi\cdot X$. Recall that $W$ defines a complex structure $I$ on $V$ such that $v\mapsto v-iI(v)$ is an isomorphism of complex vector spaces $(V,I) \to W$. Let $E$ be the real subspace of $V$ spanned by the last $g$ vectors of the standard symplectic basis. We have $E\cap I(E) =\{0\}$ since $w = I(v) \in E$ for some $v\in E$ implies $Q(w,I(w)) = Q(w,I^2(v)) = Q(v,w) = 0$ contradicting the positivity condition of a complex structure unless $v = 0$. This shows that the vectors $v-iI(v), v\in E$ span $W$ as a complex space. Thus we can find a unique basis of $W$ such that the last $g$ rows of the matrix $\Pi$ form the identity matrix. In other words, we can always assume that 
\begin{equation}\label{z}
\Pi = \begin{pmatrix}Z\\I_g\end{pmatrix},
\end{equation}
for a unique  square complex matrix $Z$ of size $g$. The conditions \eqref{riem1} and \eqref{riem2} are equivalent to the conditions
\begin{equation}\label{riem}
{}^tZ = Z, \quad \text{Im}(Z) = \frac{1}{2i}(Z-\bar{Z}) > 0.
\end{equation}
We obtain that the period space $\calD_{(2g,g)}$ parametrizing polarized AHS of weight 1 is isomorphic to the complex manifold 
\begin{equation}\label{siegel2}
 \calZ_g: = \{Z\in \Mat_g(\bbC): {}^tZ = Z, \text{Im}(Z) > 0\},
\end{equation}
called the {\it Siegel half-plane} of degree $g$. Its dimension is equal to $g(g+1)/2$. We know that $\calD_{(2g,g)}$ is a compex homogeneous space of the form $G_\bbR/K$, where 
$G_\bbR = \Aut(\bbR^{2g},Q)$  and $K$ is a compact subgroup of $G_\bbR$. Explicitly, 
$G_\bbR$ can be identified with the group of matrices (the {\it symplectic group})
$$\Sp(2g,\bbR) = \{M\in \GL(2g,\bbR): {}^tM\cdot J\cdot M = J\}.$$

The group $\Sp(2g,\bbR)$ acts on $\calZ_g$ by its natural action on isotropic subspaces $W$ of $G(g,V_\bbC)$ or on matrices $Z\in M_g(\bbC)$ by the formula
\begin{equation}\label{action}
M\cdot Z = (A\cdot Z+B)\cdot(C\cdot Z+D)^{-1},
\end{equation}
where $M$ is written as a block-matrix
$M= \tiny{\begin{pmatrix}A&B\\
C&D\end{pmatrix}}.$
Let $Z = X+iY\in \calZ_g$. Since $Y$ is a symmetric positive definite matrix, it can be written as the product $A\cdot {}^tA$ for some invertible square matrix $A$. The matrix 
$M = \tiny{\begin{pmatrix}A&X{}^tA^{-1}\\
0&{}^tA^{-1}\end{pmatrix}}$ belongs to $\Sp(2g,\bbR)$ and 
$M\cdot iI_{g} = (iA+X{}^tA^{-1})\cdot A^{t} = X+iY.$ This checks that $\Sp(2g,\bbR)$ acts transitively on $\calZ_g$.
If we take $Z$ to be the matrix $iI_g$, then the stabilizer $\Sp(2g,\bbR)_Z$ of $Z$ consists of matrices $M$ with $iA+B =i( iC+D)$, i.e. satisfying $A= D, B = -C$. The map $M\to X = A+iB$ is an isomorphism of $\Sp(2g,\bbR)_Z$ onto a subgroup of complex $g\times g$-matrices. Also the condition ${}^tM\cdot J\cdot M = J$ translates into the condition ${}^tX\cdot X = I_g$. Thus a compact subgroup $K$ could be taken to be a subgroup of $\Sp(2g,\bbR)$ isomorphic to the unitary group $U(g)$.  So we get 
\begin{equation}\label{siegel3}
\calZ_g \cong \Sp(2g,\bbR)/U(g).
\end{equation}
 A Siegel upper-half plane is an irreducible  hermitian symmetric space of type III in Cartan's classification. 

Recall that {\it hermitian symmetric  space}  is a connected complex manifold $X$ 
equipped with a hermitian metric such that  each point $x\in X$
is an isolated fixed point of a holomorphic involution  of $X$ which
preserves the metric. One can show  that the metric satisfying this
condition must be a K\"ahler metric.  The group of holomorphic isometries 
acts transitively on a hermitian symmetric space
$X$ and its connected component of the identity is a connected Lie group
$G(X)$. The isotropy subgroup of a point is a maximal compact subgroup $K$ of $G(X)$ which contains a central subgroup isomorphic to $U(1)$. An element $I$ of this  subgroup  satisfying $I^2=-1$ defines a complex structure on the tangent space of each point of $X$.  Any hermitian symmetric  space is a symmetric space
with respect to the canonical
structure of a Riemannian manifold. A  hermitian symmetric space of non-compact type is a homogeneous space $G/K$ as above with $G$ a semi-simple Lie group.  It is called irreducible if the Lie group is a simple Lie group.

The Siegel half-plane $\calZ_1$ is
 just the upper half-plane
$$\calH = \{z=x+iy\in \bbC: y > 0\}.$$
The action of $\Sp(2g,\bbR)$ on $\calZ_g$ is analogous of the Moebius transformation of the upper half-plane $z\mapsto \frac{az+b}{cd+d}$.

It is known that $\calH$ is holomorphically isomorphic to the unit 1-ball 
$$U = \{z\in \bbC: |z| < 1\}.$$
Similarly, 
the Siegel half-plane is isomorphic as a complex manifold to a bounded domain in
 $\bbC^{\frac{1}{2}g(g+1)}$ defined by 
\begin{equation}\label{bd1}\{Z\in \Mat_{g}(\bbC): {}^tZ = Z, I_g-\bar{Z}\cdot Z > 0\}\end{equation}
The isomorphism is defined by replacing matrix \eqref{z} satisfying \eqref{riem1} and \eqref{riem2} with the matrix 
$$\Pi' = \begin{pmatrix}I_g&I_g\\
iI_g&-iI_g\end{pmatrix}\cdot \Pi$$ and then reducing it back to a form \eqref{z}. 

There is a natural definition of the direct sum of abstract polarized Hodge structures which we leave to the reader. For AHS of weight 1 this of course corresponds to the operation of the direct sum of symplectic spaces and the direct sum of complex structures. It is easy to see that a polarized AHS of weight 1 decomposes into a direct sum if and only if the corresponding matrix $Z\in \calZ_g$ is the direct product of block-matrices of smaller size.

\section{Abelian varieties}\label{}

Recall that a  torus $T$ is a smooth manifold $V/\Lambda$, where $V$ is a real space of dimension $n$ and $\Lambda$ is a {\it lattice} in $V$. The additive group structure of $V$ defines a structure of a group manifold on $T$, the space $V$ acts on $T$ by translations with kernel equal to $\Lambda$. Obviously, $T$ is  diffeomorphic to $\bbR^{n}/\bbZ^{n} = (S^1)^{n}$, the product of $n$ circles.

 The space $V$ is the universal covering of $T$, and the group $\Lambda$ is identified with the fundamental group of $T$. Since it is abelian, we obtain a canonical isomorphism
$$\Lambda \cong H_1(T,\bbZ) .$$
By K\"unneth formula,
$$H_k(T,\bbZ) \cong \bigwedge^k\Lambda,\ 
H^k(T,\bbZ) \cong \bigwedge^k\Hom(\Lambda,\bbZ).$$
There is a natural bijective correspondence between tensor forms on the smooth manifold $V/\Lambda$ invariant with respect to translations and tensors on $V$. Tensoring with $\bbR$ and $\bbC$ we can identify $H^k(T,\bbR)$  with the space $k$-multilinear forms on $V= \Lambda\otimes \bbR$ and $H^k(T,\bbC)$ with the space of $k$- multilinear forms on the complexification $V_\bbC$ of $V$. 

Now assume $n = 2g$.
A complex structure on $T$ invariant with respect to translations (making $T$ a complex Lie group) is defined by a complex structure $I$ on the space $V$. The holomorphic tangent bundle of $T$ becomes isomorphic to the trivial bundle with fibre $W = (V,I)$.  A translation invariant  
 structure of a hermitian manifold on $T$ is defined by a positive definite hermitian form $h$ on $W$. It is easy to see that 
\begin{equation}\label{def2}
h(v,w)  = Q(I(v),w)+iQ(v,w),
\end{equation}
where $Q$ is a skew-symmetric form on $V$ and $g(v,w) = Q(I(v),w)$ is a positive definite symmetric form.   Since
$h(I(v),I(w)) = h(v,w)$, we see that
$g(v,w)$ is a Riemannian metric on $T$ invariant with respect to the complex structure and also $I$ is an isometry of the symplectic space $(V,Q)$. It follows from \eqref{def2} that the conjugate complex structure $-I$ on $V$ is positive with respect to $Q$. Thus a translation invariant structure of a hermitian manifold on $T$ defines a Hodge structure of weight 1 on $V$ with polarization $Q$ and $W=V^{10}$. As is easy to see the converse is  true. 

 Extend $h$ to a hermitian form $H$ on $V_\bbC$ by requiring that the subspaces $W$ and $\overline{W}$ are orthogonal with respect to $H$ and 
$H|W = h,\  H(\bar{w},\bar{w}') = -H(w,w') $. Write $w = \frac{1}{2}(v-iI(v)), w' = \frac{1}{2}(v'-iI(v'))$ for some $v,v'\in V$. Since
 $$2iQ(w,\bar{w}') = 2iQ(\frac{1}{2}(v-iI(v)),\frac{1}{2}(v'+iI(v')) = Q(v,I(v'))+iQ(v,v') = h(v,v'),$$ we see that 
$H$ is defined by
$$H(x,y)  = 2iQ(x,\bar{y})$$
and hence $H$ and $\Omega = -2Q$ are associated hermitian and skew-symmetric forms on $V_\bbC$ as is defined in the previous section. Let $z_\alpha$ be the coordinate functions on $W$ with respect to some basis $e_1,\ldots,e_g$ and $\bar{z}_\alpha$ be the coordinate functions on $\overline{W}$ with respect to the conjugate basis $\bar{e}_1,\ldots,\bar{e}_g$. Then 
$z_\alpha,\bar{z}_\alpha$ are coordinate functions on $V_\bbC$ satisfying
$$\overline{z_\alpha(x)} = \bar{z}_\alpha(\bar{x}).$$ 
Let $h_{\alpha\beta} = H(e_\alpha,e_\beta)$ so that, for any $x,y\in W$, 
$$h(x,y) = \sum h_{\alpha\beta}z_\alpha(x)\bar{z}_\beta(y).$$
Then $\Omega(e_\alpha,e_\beta) = \Omega(\bar{e}_\alpha,\bar{e}_\beta) = 0$ and
$$\Omega(e_\alpha,\bar{e}_\beta) = ih(e_\alpha,e_\beta) = ih_{\alpha\beta}.$$
This shows that, for any $x,y\in V_\bbC$,
$$\Omega(x,y) = i\sum h_{\alpha\beta}z_\alpha(x)\wedge \bar{z}_\beta(y).$$
Comparing this with \eqref{metric}, we see that $\Omega$ defines the K\"ahler form on $T$ associated to the hermitian metric $h$. Since it has constant coefficients, its closedness is obvious. 

Assume that  the  symplectic form $Q$ satisfies the integrality condition with respect to the lattice $\Lambda$, then the K\"ahler form $\omega$ associated to the hermitian metric is a Hodge form, and hence  the complex torus $X = W/\Lambda$ admits an embedding in a projective space such that the cohomology class in $H^2(T,\bbZ)$ which is dual (with respect to the Poincar\'e duality) to the homology class in $H_{2n-2}(T,\bbZ)$ of a hyperplane in $\bbP^n$ is equal to some positive multiple of the cohomology class $[\omega]$.

Conversely, let $T= W/\Lambda$ be a complex torus. We identify the real vector space $V$ underlying $W$ with  $\Lambda\otimes \bbR$ and $V_\bbC$ with $\Lambda\otimes \bbC$. A choice of  a symplectic  form $Q$ on $V$  such that the complex structure $W$ on $V$ is positive with respect to $-Q$ is called a {\it polarization} on $T$. Two polarizations are called equivalent if the symplectic forms differ by a constant factor equal to a positive rational number. A polarization defines a hermitian form on $W$ and a structure of a hermitian complex manifold on $T$. A pair $(T,Q)$ is called a {\it polarized compex torus}. If additionally $Q(\Lambda\times \Lambda)\subset \bbZ$, then the polarized torus is called a {\it polarized abelian variety}. Recall that an {abelian variety} (over compex numbers) is a projective algebraic variety  isomorphic, as a complex manifold, to a complex torus. A projective embedding of an abelian variety $X$ defines a structure of a polarized abelian variety on $X$, the Hodge class being the cohomology class of a hyperplane section.

Not every complex torus $W/\Lambda$  is an abelian variety. A polarization $Q$ is defined by an integral matrix $A$ in a basis of $\Lambda$. If we identify $V_\bbC$ with $\bbC^{2g}$ by means of this basis and the space $\bar{W}$ with the column space of a matrix $\Pi$ (called the {\it coperiod matrix} of the torus), then the necessary and sufficient conditions of positivity of the complex structure of $\bar{W}$ with respect to $Q$ are the conditions
\begin{equation}\label{riem}
{}^t\Pi\cdot A\cdot \Pi = 0, \quad -i{}^t\Pi\cdot A\cdot \bar{\Pi} > 0.
\end{equation}

These conditions can be combined by introducing the square matrix $\tilde{\Pi} = (\Pi|\bar{\Pi})$ (corresponding to the matrix $\Pi_0$ in (PM1)-(PM3)). The condition is 
\begin{equation}\label{bigriem}
-i{}^t\tilde{\Pi}\cdot A\cdot \overline{\tilde{\Pi}} = \begin{pmatrix}M&0_g\\
0_g&-{}^tM\end{pmatrix},
\end{equation}
where $M$ is a positive definite hermitian matrix. 

Since $V_\bbC = W\oplus \overline{W}$, the matrix  $\tilde{\Pi}^{-1}$ is invertible. Write 
$$\tilde{P} = \tilde{\Pi}^{-1} = (P|\bar{P})$$
for some $2g\times g$-matrix. The $j$th column of $P$ gives an expression of the $j$th basis vector of $\Lambda$ as a linear combination of a basis of $W$ formed by the columns of $\Pi$.  It is called the {\it period matrix} of $T$. One can restate the condition \eqref{bigriem} in terms of the period matrix by
\begin{equation}\label{bigriem'}
i{}^t\overline{\tilde{P}}\cdot A^{-1}\cdot {}^t\tilde{P}  = \begin{pmatrix}M'&0_g\\
0_g&-{}^tM'\end{pmatrix},
\end{equation}
where $M' = M^{-1}$ is a positive definite hermitian matrix.  These  are so called the {\it Riemann-Frobenius conditions}. 

\begin{ex} Any one-dimensional torus is an abelian variety (an elliptic curve). In fact let 
$T$ be isomorphic to the complex torus $\bbC/\bbZ\gamma_1+\bbZ\gamma_2$. Under the multiplication by $\gamma_1$, the torus is isomorphic to the torus $\bbC/1\cdot \bbZ+\tau\bbZ$, where $\tau = \gamma_2/\gamma_1$. Replacing $\tau$ by $-\tau$, if necessary we may assume that $\text{Im}(\tau) > 0$. Thus the period matrix of $T$ with respect to the basis $1$ of $\bbC$ and the basis $1,\tau$ of $\Lambda$ is equal to
$P = (1,\tau)$. The Riemann-Frobenius condition is the existence of an integer $r$ such that  
$$i\begin{pmatrix}1&\bar{\tau}\\
1&\tau\end{pmatrix}\cdot \begin{pmatrix}0&1/r\\
-1/r&0\end{pmatrix}\cdot \begin{pmatrix}1&1\\
\tau&\bar{\tau}\end{pmatrix} = \begin{pmatrix}a&0\\
0&-a\end{pmatrix},$$
where $a$ is a positive real number. Computing the product we obtain that the condition is equivalent to the condition 
$ir^{-1}(\bar{\tau}-\tau))\in \bbR_{> 0}$ which is satisfied for any negative integer $r$.  

We leave to the reader to check that the complex torus $T = \bbC^2/P\cdot \bbZ^4,$ where
$$P = \begin{pmatrix}\sqrt{-2}&\sqrt{-3}\\
\sqrt{-5}&\sqrt{-7}\end{pmatrix}$$
does not satisfy condition \eqref{bigriem'} for any integral skew-symmetric matrix $A$. Thus it is not a projective algebraic variety.
\end{ex}

Let $(W/\Lambda,Q)$ be a polarized abelian variety.  It is known from linear algebra that $\Lambda$ admits a basis such that the matrix of $Q|\Lambda$ is equal to 
\begin{equation}\label{matrix}
J_D = \begin{pmatrix}0_g&D\\
-D&0_g\end{pmatrix},\end{equation}
where $D = \text{diag}[d_1,\ldots,d_g]$ with $0 < d_1|d_2|\ldots |d_g$. The vector $(d_1,\ldots,d_g)$  is an invariant of an integral valued  skew-symmetric form on $\Lambda$. The equivalence class of $D$ with respect to the equivalence relation defined by $D\sim D'$ if $D' = aD$ for some positive rational number $a$  is called the {\it type} of the polarization of $(W/\Lambda,Q)$.  We will always represent it by a primitive $D$ (i.e. not positive multiple of any other $D$). A polarization of type  $(1,\ldots,1)$ is called a {\it principal polarization}.

Two polarized abelian varieties $(W/\Lambda,Q)$ and $(W'/\Lambda',Q')$ are called isomorphic if there exists an isomorphism of complex spaces $f:W\to W'$ such that $f(\Lambda) = \Lambda'$ and $Q'\circ (f\times f):\Lambda\times \Lambda\to \bbZ$ is equivalent to $Q$. Clearly, isomorphic polarized varieties have the same type of polarization. Any polarized abelian variety of type $(d_1,\ldots,d_g)$ is isomorphic to the polarized variety $\bbC^g/P\cdot \bbZ^{2g}$, where $P$ is the period matrix such that the matrix of $Q$ with respect to the basis of $\Lambda$ formed by the columns of $P$ is the matrix \eqref{matrix}. Then $P$ satisfies \eqref{bigriem'}. The matrix
$$\Pi' = \begin{pmatrix}1_g&0_g\\
0_g&D\end{pmatrix}\cdot {}^tP$$
satisfies \eqref{riem1} and \eqref{riem2}.
 As we have seen, there exists an invertible complex matrix $X$ such that $\Pi\cdot X$ is of the form \eqref{z}, where $Z\in \calZ_g$. 
This implies that our torus is isomorphic to the torus $\bbC^g/P\cdot \bbZ^{2g}$, where $P = (Z|D)$ and $Z\in \calZ_g$.
The period matrix of this form is called a {\it normalized period matrix}.  Two polarized abelian varieties with normalized period matrices $P = (Z|D)$ and $P' = (Z'|D)$ are isomorphic if and only if there exists a complex matrix $X$ defining a map $\bbC^g\to \bbC^g$ such that 
$$X\cdot (Z|D) = (Z'|D)\cdot M,$$
where $M$ is an invertible integral $2g\times 2g$-matrix satisfying
$${}^tM\cdot J_D\cdot M = J_D.$$
Let $\Sp(2g,\bbZ)_D \subset \GL(2g,\bbZ)$ be the group of such matrices. If $D = I_g$, i.e., $J_D = J$, the group $\Sp(2g,\bbZ)_D$ is the {\it Siegel modular group} 
$$\Gamma_g = \Sp(2g,\bbZ)$$
of symplectic integral matrices.
For any $M\in \Sp(2g,\bbZ)_D$, the conjugate matrix 
$$N = \begin{pmatrix}1_g&0_g\\
0_g&D\end{pmatrix}\cdot M\cdot \begin{pmatrix}1_g&0_g\\
0_g&D\end{pmatrix}^{-1}$$
belongs to $\Sp(2g,\bbQ)$ and leaves invariant the lattice $\Lambda_D$ in $\bbR^{2g}$ with period matrix $ \begin{pmatrix}1_g&0_g\\
0_g&D\end{pmatrix}$. It is easy to see that any such matrix arises in this way from a matrix from $\Sp(2g,\bbZ)_D$. Denote the group of such matrices by $\Gamma_D$. We see that 
$$\Gamma_D = \begin{pmatrix}1_g&0_g\\
0_g&D\end{pmatrix}\cdot \Sp(2g,\bbZ)_D\cdot \begin{pmatrix}1_g&0_g\\
0_g&D\end{pmatrix}^{-1}.$$

Assume that 
$P = (Z|D),P' = (Z'|D)$ are normalized. The corresponding tori are isomorphic as polarized abelian varieties if and only if there exists a matrix $M\in \Sp(2g,\bbZ)_D$ and a matrix 
$X\in \GL(g,\bbC)$ such that
\begin{equation}\label{eqq}
X\cdot (Z'|D) = (Z|D)\cdot M.\end{equation}
Write $M$ in the form
\begin{equation}\label{conj}
M = \begin{pmatrix}1_g&0_g\\
0_g&D^{-1}\end{pmatrix} \cdot N\cdot \begin{pmatrix}1_g&0_g\\
0_g&D\end{pmatrix},
\end{equation}
where 
$$N = \begin{pmatrix}N_1&N_2\\
N_3&N_4\end{pmatrix}\in \Gamma_g(D),$$
and the blocks are all of size $g$.
Then \eqref{eqq} is equivalent
$$X = Z\cdot N_3+N_4, \ X\cdot Z' = Z\cdot N_1+N_2.$$
This shows that $X$ is determined by $N$ and 
\begin{equation}\label{act1}
Z' = (Z\cdot N_3+N_4)^{-1}\cdot (Z\cdot N_1+N_2).\end{equation}
Transposing these matrices  and using that ${}^tZ = Z$ because $Z\in \calZ_g$, we obtain
\begin{equation}\label{act2}
Z' = ({}^tN_1\cdot Z+{}^tN_2)\cdot ({}^tN_3\cdot Z+{}^tN_4)^{-1}.\end{equation}

The group $\Sp(2g,\bbQ)$ is invariant with respect to the transpose operation, let 
${}^t\Gamma_D$ be its subgroup of matrices ${}^tM$, where $M\in \Gamma_D$. Then 
$\Gamma_D$ acts on $\calZ_g$ by the formula \eqref{act1} or \eqref{act2} and ${}^t\Gamma_D$ acts on $\calZ_g$ by the restriction of the action \eqref{action} of $\Sp(2g,\bbR)$ on $\calZ_g$.

Summarizing, we obtain the following.

\begin{thm} The natural map assigning to a polarized abelian variety its normalized period matrix defines a bijection between isomorphism classes of abelian varieties with polarization of type $D = (d_1,\ldots,d_g)$ and the orbit space
$$\calA_D = \Gamma_D\backslash \calZ_g.$$

\end{thm}

In fact one can show more, the coarse moduli space of abelian varieties with polarization of type $D$ exists and is isomorphic to $\Gamma_D\backslash \calZ_g$.

\section{Picard and Albanese varieties}

Let $X$ be a compact K\"ahler manifold.   It defines a Hodge structure of weight 1 
$$V_\bbC = H^1(X,\bbC) = H^{10}(X)\oplus H^{01}(X)$$
on the real space $V = H^1(X,\bbR)$ of dimension $b_1(X) = 2q$. The map
$$\alpha:H_1(X,\bbZ)\to H^{10}(X)^*, \quad \alpha(\gamma)(\omega) = \int_\gamma\omega,$$
defines an isomorphism from $H_1(X,\bbZ)/\Tors$ onto a lattice $\Lambda$ in $H^{10}(X)^* = H^0(X,\Omega_X^1)^*$. The complex torus 
$$\Alb(X) = H^{10}(X)^*/\Lambda$$
is called the {\it Albanese torus} of $X$. Its period matrix with respect to a basis $\omega_1,\ldots,\omega_q$ in $H^{10}(X)$ and a basis $(\alpha(\gamma_1),\ldots,\alpha(\gamma_{2q}))$ of $\Lambda$ is the matrix
\begin{equation}\label{per}
P = (\int_{\gamma_j}\omega_i).
\end{equation}
Fixing a point $x_0\in X$ defines a natural map (the {\it Albanese map})
$$\alb_{x_0}:X\to \Alb(X), \quad x \to (\int_{x_0}^x\omega_1,\ldots,\int_{x_0}^x\omega_g)\  \mod \Lambda,$$
where the integration is taken with respect to any 1-chain originated from $x_0$ and ending at $x$. Since the difference of two such chains is a 1-cycle, the point in the torus does not depend on a choice of a path. 
The Albanese map satisfies the following universal property: for any holomorphic map $\phi:X\to T$ to a complex torus $T$ there exists a unique (up to translation) map of complex tori $f:\Alb(X)\to T$ such that $f\circ \alb = \phi$.

One can define another complex torus associated to $X$. Consider the projection $\Lambda'$ of  $H^1(X,\bbZ)\subset V$ to $H^{01}(X)$. It is a lattice in the complex space $H^{01}(X)$. The complex torus
\begin{equation}
\Pic^0(X) = H^{01}(X)/\Lambda'
\end{equation}
is called the {\it Picard torus} of $X$. Recall that $H^{01}(X) \cong H^1(X,\calO_X)$, where $\calO_X$ is the sheaf of holomorphic functions on $X$. The exponential exact sequence
$$0\to \bbZ \to \bbC\overset{\tiny{z\mapsto e^{2\pi iz}}}{\longrightarrow}\bbC^* \to 0$$
defines an exact sequence of sheaves 
$$0\to \bbZ \to \calO_X\overset{\exp}{\longrightarrow} \calO_X^* \to 0$$
and the corresponding exact cohomology sequence
$$H^1(X,\bbZ) \to H^1(X,\calO) \to H^1(X,  \calO_X^*)  \overset{c}{\to} H^2(X,\bbZ).$$
One can show that the first map  coincides with the composition of the maps  
$H^1(X,\bbZ) \to H^1(X,\bbC)\to H^{01}(X)$ so that 
$$\Pic^0(X) \cong \Ker(H^1(X,\calO_X^*)\overset{c}{\to} H^2(X,\bbZ)).$$
The group $H^1(X,\calO_X^*)$ is the {\it Picard group} of $X$ of isomorphism classes of holomorphic line bundles over $X$. The coboundary map $c$ is interpreted as taking the first Chern class of a line bundle.  Thus the Picard variety is a complex torus whose set of points is naturally bijective to the set of isomorphism classes of holomorphic line bundles $L$ with trivial first Chern class $c_1(L)$.

Of course, one can also consider the complex tori $\Alb(X)'$ and $\Pic^0(X)'$, interchanging  $H^{10}(X)$ with $H^{01}$. This  just replaces the complex structure of the torus  to the conjugate one.  

The complex tori $\Alb(X)$ and $\Pic^0(X)$ are examples of {\it dual complex tori} in the sense of the following definition.  Let $T = V/\Lambda$ be a torus and $W = (V,I)$ be a complex structure on $V$. Recall that we identify $W$ with the complex subspace of $V_\bbC$ of vectors 
$v-iI(v)$. For any $\bbR$-linear function $\phi\in V^*$ and $w= v+iI(v)\in \overline{W}$ set 
$$\tilde{\phi}(w) = \phi(I(v))+i\phi(v).$$
It is immediately checked that $\tilde{\phi}$  is a $\bbC$-linear function on $\overline{W}$ with imaginary part equal to $\phi$ (when we identify $\overline{W}$ with $(V,-I)$). In this way we identify the real spaces $V^*$ and $\overline{W}^*$ and hence consider $\overline{W}$ as a complex structure on $V^*$.  Let 
$$\Lambda^* = \{f\in \overline{W}^*:\text{Im}(f)(\Lambda)\subset \bbZ\}.$$
The dual complex torus is
$$T^*: = \overline{W}^*/\Lambda^*.$$ 
In our case, $H^{10}(X)^* = \overline{H^{01}(X)}^*.$ Since 
$H^1(X,\bbZ) = \{\phi\in H^1(X,\bbR):\phi(H_1(X,\bbZ))\subset \bbZ\}$, we see that  $\Lambda' = \Lambda^*$.

The dual complex tori correspond to the dual Hodge structure on $V^*$
$$V_\bbC^* = (V^*)^{10}\oplus (V^*)^{01},$$
where  
$$(V^*)^{10} = (V^{01})^* = (V_\bbC/V^{10})^* = (V^{10})^\perp,$$ 
$$(V^*)^{01} = (V^{10})^* = (V_\bbC/V^{01})^* = (V^{01})^\perp.$$ 

Now let $Q$ be a polarization of the Hodge structure on $V$  with hermitian form $H$. It  defines a dual symplectic form $Q^*$ on $V^*$ by viewing $Q$ as a bijective map $V\to V^*$ and setting $Q^* = -Q^{-1}:V^*\to V$. Since the subspaces $V^{10},V^{01}$ are isotropic with respect to $Q$, $Q(V^{10}) = (V^{10})^\perp = (V^{01})^*$, 
$Q(V^{01}) = (V^{10})^*$. This shows that $(V^{10})^*$ and $(V^{01})^*$ are isotropic with respect to $Q^*$.  Let $V^{10}$ defines a complex structure $I$ on $V$. Then the complex structure on $(V^{10})^*$ is defined by the operator $I^*$ which is adjoint to $I$ with respect to the natural pairing between $V$ and $V^*$. The  symmetric bilinear form $Q(v,I(v'))$ can be considered as the composition of the map $V\to V^*$ and $I^*:V^*\to V^*$. The symmetric form $Q^*(I^*(\alpha),\beta)$ defines the inverse map. 
This shows that $I$ is positive with respect to $Q$ if and only if $Q^*$ is positive with respect to $I^*$. Thus $Q^*$ defines a polarization of the dual of a polarized AHS of weight 1. It is called the {\it dual polarization}. 

An example of dual Hodge structures is a Hodge structure on odd-dimensional cohomology  $H^k(X,\bbR)$ and on $H^{2n-k}(X,\bbR)$, where $X$ is a K\"ahler manifold of dimension $n$. The duality is defined by the 
Poincar\'e duality. It is also a duality of polarized Hodge structures.

It is clear that a polarized Hodge structure on $H^1(X,\bbR)$ defines a polarization on $\Pic^0(X)$ and the dual polarization on $\Alb(X)$. An integral polarization defined by a Hodge class makes these tori polarized abelian varieties.

Suppose $\Lambda \subset V$ is an integral structure of $Q$. Choose a basis $(v_i)$ in $\Lambda$ such that the matrix of $Q$ is equal to a matrix $J_D$ for some $D$.  Let $(v_i^*)$ be the dual basis in $V^*$. The map $Q:V\to V^*$ sends $v_i$ to $d_iv_{i+g}^*$  and $v_{i+g}$ to $-d_iv_{i}^*$ for $i = 1,\ldots,g$. Thus the matrix of $Q^*$ with respect to the basis 
$v_{1+g}^*,\ldots,v_{2g}^*, v_1^*,\ldots,v_g^*$ is equal to $J_D^{-1}$. Multiplying $Q^*$ by $d_g$ we get an integral structure of the dual AHS of type $D^* = (1,d_g/d_{g-1},\ldots,d_g/d_1)$ with respect to the dual lattice $\Lambda^*$. 

In particular, we see that the dual of a polarized  abelian variety of type $D$ is a polarized abelian variety of type $D^*$.

Now let $Q$ be a polarization on $T$ with the integral structure defined by $\Lambda$. Let $H$ be the associated hermitian form on $V_\bbC$ and let $h$ be its restriction to $W$. It defines a $\bbC$-linear map 
$h:W\to \overline{W}^*, x \mapsto H(x, )$. Since $\text{Im}(h)(x,y) = Q(x,y)$, we see that 
$h(\Lambda) \subset \Lambda^*$. Thus we have a holomorphic map of complex tori
$$h:T \to T^* = \Pic^0(T).$$ 
Obviously the map is surjective, and its kernel is isomorphic to the finite abelian group of order $\det J_D = (d_1\cdots d_g)^2$ isomorphic to the cokernel of the map $Q:\Lambda\to \Lambda^*$. In particular, $h$ is an isomorphism if $D = 1_g$.

Let $P$ be the period matrix of $T$ with respect to the chosen bases of $\Lambda$ and $W$. It defines a map $V_\bbC = \Lambda_\bbC \to W$, the direct sum of the transpose map $W^*\to V_C^*$ and its  conjugate map $\overline{W}^* \to V_C^*$ define the dual Hodge structure. Thus the matrix ${}^t\bar{P}$ is the coperiod matrix $\Pi^*$ of the dual torus $T^*$ with respect to the bases 
$(v_i^*)$ and $(\bar{w}_i)$. The period matrix $P^*$ of $T^*$ is reconstructed from the equation
$$\begin{pmatrix}P^*\\
\bar{P}^*\end{pmatrix} = (\Pi^*|\bar{\Pi}^*)^{-1} = ({}^t\bar{P}|{}^tP)^{-1}.$$

\begin{ex} Let $T = W/\Lambda$ be a complex torus. The space $T\times W$ is identified with the trivial holomorphic tangent bundle of $T$. Thus 
the dual space $W^*$ is the space of holomorphic differentials $H^{10}(X)$. Thus 
$$T \cong  \Alb(T).$$
Hence 
$$\Pic^0(T) = T^*.$$
Also we see that a polarization on $T$ defines the dual polarization of $\Pic^0(T)$. In particular, if $T$ is an abelian variety with a principal polarization, then $T \cong \Pic^0(T)$ as polarized abelian varieties.
\end{ex}

So far we have defined complex tori associated to a polarized  AHS of weight $1$ with integral structure $\Lambda \subset V$. One can associate complex tori to an AHS of any odd weight $k = 2m+1$. In fact, set 
\begin{equation}\label{hdg}
V^{10} =  \bigoplus_{s=0}^{m} H^{m+1+s,m-s}(X), \quad V^{01} =  \bigoplus_{s=0}^{m} H^{s,k-s}(X),
\end{equation}
then $V_\bbC = V^{10}\oplus V^{01}$ is an AHS of weight 1 which is polarized with respect to the polarization form $Q$ of the original AHS. This defines two dual complex tori 
$$J = V^{10}/\Lambda,\  \check{J} = (V^{01})^*/\Lambda^* $$
where $\Lambda$ is the projection  of $\Lambda$ to $V^{10}$. The torus $J$ is called the {\it Griffiths complex torus} associated to AHS. 

In the case when  the Hodge structure is the Hodge structure on odd-dimensional cohomology  $H^{k}(X,\bbR)$ of a $n$-dimensional K\"ahler variety with integral structure defined by $\Lambda = H^{k}(X,\bbZ)$ we obtain the definition of the {\it kth Griffiths intermediate Jacobian} $J = J^{k}(X)$. We have
$$\check{J}^{k}(X) \cong J^{2n-k}(X).$$
 In particular, they coincide when $k = n = 2m+1$. In this case $J^n(X)$ is called the {\it intermediate Jacobian} of $X$. Clearly,
$$J^1(X) = \Pic^0(X),\quad J^{n-1}(X) = \Alb(X).$$

Note that, in general,  $Q$ does not define a polarization of the AHS \eqref{hdg} since the hermitian form $iQ(x,\bar{x})$ is not positive on $V^{10}$; it changes signs on the direct summands. Thus the Griffiths tori are not abelian varieties in general. However, in one special case they are. This is the case when 
all summands except one in $V^{10}$ are equal to zero. We call such AHS an {\it exceptional AHS}. In this case $J$ is a polarized abelian variety. It is a principally polarized variety if $Q$ is unimodular on $\Lambda$. For example, the  intermediate Jacobian $J^{n}(X)$ of an odd-dimensional K\"ahler manifold with exceptional Hodge structure on the middle cohomology is a principally polarized abelian variety. For example, any K\"ahler 3-fold  with $h^{30} = 0$ (it must be an algebraic variety in this case) defines a principally polarized abelian variety 
$J^3(X)$ of dimension $h^{21}$. Another example  is a  rigid Calabi-Yau 3-fold for which $h^{30} = 1, h^{21} = 0$. 

Let $T$ be a connected complex variety. A {\it family of polarized abelian varieties} with base $T$ is a closed  subvariety $Y$ of  $\bbP^N\times T$ such that the each fibre $X_t$ of the second projection $f:Y\to T$ is an abelian variety of  dimension $g$. Under the natural map $X_t \to \bbP^N\times\{t\}$, each fibre is isomorphic to a closed projective subvariety of $\bbP^N$ and hence acquires a structure of a polarized abelian variety of some type $D$ independent of $t$. The period map is defined and is a holomorphic map
$$\bar{\phi}:T\to \Gamma_D\backslash \calZ_g.$$
If we identify the orbit space with the coarse moduli space $\calA_g(D)$ of abelian varieties with polarization of type $D$, then the map is $X_t\mapsto \Pic^0(X)$. The abelian variety $X$ can be uniquely reconstructed from $\Pic^0(X)$ as the dual polarized abelian variety. This proves the Global Torelli Theorem for polarized abelian varieties.  It gives an isomorphism
$$\per:\calA_g(D) \to \Gamma_{D^*}\backslash \calZ_g \cong \calA_{D^*}.$$

\begin{ex} Let $M$ be a compact smooth oriented  2-manifold. Then $b_1(M) = 2g$, where $g$ is the genus of $M$. We choose a complex structure on $M$ which makes it a compact Riemann surface $X$, or a projective nonsingular curve.   It is known that the group $H_1(M,\bbZ)$ is a free abelian group of rank $2g$ which admits a standard symplectic basis $\gamma_1,\ldots,\gamma_{2g}$ with respect to the cup-product 
$$Q:H^1(M,\bbZ)\times H^1(M,\bbZ) \to H^2(M,\bbZ) \cong \bbZ.$$
This defines an integral  polarization of the Hodge structure corresponding the the K\"ahler class generating $H^2(X,\bbZ) =\bbZ$. It gives a principal polarization on $\Pic^0(X)$. Hence $\Alb(X) \cong \Pic(X)$ are isomorphic as principally polarized abelian varieties and coincide, by definition, with the {\it Jacobian variety} $\Jac(X)$ of $X$. We know that there exists a basis $\omega_1,\ldots,\omega_g$ in the space  of holomorphic 1-forms $H^{10}(X)$ such that the period matrix of $\Jac(X)$
$$\Pi = (\pi_{ij}), \quad \pi_{ij} = \int_{\gamma_i}\omega_j$$
is normalized, i.e. has the form ${}^t(Z|I_g)$, where $Z\in \calZ_g$. 

A {\it marked Riemann surface} is a Riemann surface $X$ together with a choice of a symplectic basis in $H^1(X,\bbZ)$. An isomorphism of marked Riemann surfaces $(X,(\gamma_i)) \to (X',(\gamma_i'))$ is an isomorphism $f:X\to X'$ of complex varieties such that $f^*(\gamma_i') = \gamma_i$. The {\it Torelli Theorem} for Riemann surfaces asserts that  
two marked Riemann surfaces are isomorphic if and only if their normalized period matrices with respect to the symplectic matrices are equal. In other words they define the same point in the Siegel half-plane $\calZ_g$. Two different markings on the same surface define the points in $\calZ_g$ equivalent with respect to the action of the Siegel modular group $\Gamma_g = \Sp(2g,\bbZ)$. Thus we see that two Riemann surfaces are isomorphic if and only if the corresponding Jacobian varieties are isomorphic as principal abelian varieties, i.e. define the same point in $\calA_g = \calZ_g/\Gamma_g$. The coarse moduli space $\calM_g$ of Riemann surfaces of genus $g$ exists and is an algebraic variety of dimension $3g-3$. The period map
$$\per:\calM_g \to \Gamma_g\backslash \calZ_g \cong \calA_g$$
is just the map $X\mapsto \Jac(X)$. Thus the Torelli Theorem for Riemann surfaces is the statement that two Riemann surfaces are isomorphic if and only if their Jacobian varieties are isomorphic as principally polarized abelian varieties. Clearly, this also establishes the Global Torelli Theorem for Hodge structures on one-dimensional cohomology of Riemann surfaces.

In case $g = 1$, we have an isomorphism $\calA_1 = \calH/\SL(2,\bbZ) \cong \bbC$. It is given by the {\it absolute invariant} function $j(\tau)$. The Weierstrass function
$$\wp(z) = z^{-2}+\sum_{(m,n)\in \bbZ^2\setminus\{0\}}\bigl(\frac{1}{(z+m+n\tau)^2}-\frac{1}{(m+n\tau)^2}\bigr)$$
defines a holomorphic map 
$$\bbC\setminus (\bbZ+\tau\bbZ)\to \bbC^2, \quad \bbZ \mapsto (\wp(z),\wp'(z))$$
which can be extended to an isomorphism from the elliptic curve $E_\tau = \bbC/\bbZ+\tau\bbZ$ to  the plane cubic curve in $\bbP^2$ defined by
$$x_2^2x_0-4x_1^3+g_2x_1x_0^2+g_3x_0^3 = 0.$$
We have 
$$j(\tau) = 1728\frac{g_2^3}{4g_2^3+27g_3^2}$$ and
$E_\tau\cong E_{\tau'}$ if and only if $j(\tau) = j(\tau').$
\end{ex}

\section{Eigenperiods}\label{sec6}
Let $V_\bbC = \oplus V^{pq}$ be a polarized AHS of weight $k$ on a vector space $(V,Q)$. A {\it Hodge isometry} of $V_\bbC$ is a $\bbR$-linear  automorphism of $(V,Q)$ such that its linear extention to $V_\bbC$ preserves the Hodge decomposition. Let $A$ be a finite abelian group acting on $(V,Q)$ by Hodge isometries via a linear representation $\rho:A\to \Aut(V,Q)$.  The vector 
space $V_\bbC$ splits into a direct sum of eigensubspaces 
$$V_\bbC(\chi) = \{x\in V_\bbC:g(x) = \chi(g)x, \forall g\in A\},
$$
corresponding to different characters $\chi\in \hat{A}:=\Hom(A,\bbC^*)$. 
The reality condition of the representation $\rho$ on $V_\bbC$ is equivalent to the condition
\begin{equation}\label{real}
\overline{V_\chi} =  V_{\bar{\chi}}.
\end{equation} 
We have a decomposition
\begin{equation}\label{dec1}
V_\bbC(\chi) = \bigoplus_{p+q=k} V^{pq}_\chi, \quad V^{pq}_\chi = V_\bbC(\chi)\cap V^{pq}.
\end{equation}
It is not a Hodge decomposition in general because
$$\overline{V^{pq}_\chi} = V^{qp}_{\bar{\chi}} $$
and also because the complex vector space $V_\bbC(\chi)$ does not have a natural identification with a complexification of a real vector space. 
However it satisfies properties (HD2) and (HD3) of a polarized AHS of weight $k$.  
Since, for any $x\in V_\bbC(\chi), y\in V_\bbC(\chi')$,
\begin{equation}\label{orth}
Q(x,y) = Q(g(x),g(y)) = Q(\chi(g)x,\chi'(g)y) = \chi(g)\chi'(g)Q(x,y),
\end{equation}
we obtain that the  eigensubspaces $V_\bbC(\chi)$ and  $V_\bbC(\chi')$ are orthogonal if  $\chi \ne \overline{\chi}'.$ 

Denote by $(F_\chi^p)$ the flag of subspaces 
$$F_\chi^p = \oplus_{p'\ge p} V_\chi^{p'q}.$$
It satisfies the properties
\begin{itemize}
\item[(i)] $F_\chi^p\cap \overline{F_\chi^{k-p+1}} = \{0\}$;
\item[(ii)] $(-1)^{\frac{k(k-1)}{2}}Q(Cx,\bar{x}) > 0$, where $C$ acts on $V_\chi^{pq}$ as multiplication by $i^{p-q}$.
\end{itemize}
Note that  similar to \eqref{rec} the spaces $V^{pq}_\chi$ can be reconstructed from the flag $(F^p_\chi)$ 
\begin{equation}\label{rec'} 
V_\chi^{pq} = \{x\in F^p_\chi:Q(x,\bar{y}) = 0, \ \forall y\in F_\chi^{p+1}\}.
\end{equation}
Let 
$$f_p(\chi)= \dim F_\chi^p, \quad \bff(\chi) =  (f_0(\chi),\ldots,f_k(\chi)).$$
Let $H(x,y) = Q(x,\bar{y})$ be the hermitian form on $V_\bbC(\chi)$.  We define $\Fl(\bff(\chi),V_\bbC(\chi))^0$ to be the open subset of the flag variety $\Fl(\bff(\chi),V_\bbC(\chi))$ which consists of flags $(F^p)$  satisfying conditions (i) and (ii).

If $\chi$ is real, i.e. $\chi(g) = \pm 1, $, we define $\Fl(\bff(\chi),V_\bbC(\chi))^0$ to be the period space 
$\calD_{\bff(\chi)}(V(\chi),Q|V(\chi)).$ Note that in this case the decomposition \eqref{dec1} is a polarized Hodge structure on the real eigenspace $V(\chi)$. 

Let $\wt(\rho)$ denote the subset of $\chi\in \hat{A}$ such that $(V_\bbC)(\chi) \ne 0$. 
Let $\wt_0(\rho)\subset \hat{A}$ be the subset of real characters.  Then 
$\wt(\rho)$ can be written as a disjoint union of three subsets 
\begin{equation}\label{partition}
\wt(\rho) = \wt_+(\rho)\coprod \wt_-(\rho)\coprod \wt_0(\rho) ,
\end{equation}
such that $\overline{\wt_\pm(\rho)} = \wt_\mp(\rho)$. The conjugation map of $V_\bbC$   sends the set $\Fl(\bff(\chi),V_\bbC(\chi))$ to the set
$\Fl(\bff(\bar{\chi}),V_\bbC(\bar{\chi}))$.
The group $A$ acts naturally on the period space $\calD_\bff(V,Q)$ via its representation $\rho$. Let
$$\calD_\bff(V,Q)^\rho = \{x\in \calD_\bff(V,Q): \rho(a) (x) = x, \forall a\in A\}.$$
We have a natural map
\begin{equation}\label{product}
\Phi:\calD_\bff(V,Q)^\rho \longrightarrow \prod_{\chi\in \wt_+(\rho)\cup \wt_0(\rho)}\Fl(\bff(\chi),V_\bbC(\chi))^0.
\end{equation}
It is obviously injective since one can reconstruct the Hodge structure from the image using \eqref{rec'} applied to $F_\chi^p$ and $\overline{F_{\bar{\chi}}^p}$. 

For any  K\"ahler manifold $X$ together with a K\"ahler class $[\omega]$ let $\Aut(X,[\omega])$ denote the group of holomorphic automorphisms leaving $[\omega]$ invariant. Fix a real vector space $V$ together with a non-degenerate bilinear form $Q_0$, symmetric (resp. skew-symmetric) if $k$ is even (resp. odd) and a linear representation $\rho_0:A\to \GL(V,Q_0)$. A {\it $\rho$-marking} of a polarized K\"ahler manifold $X$ on which $A$ acts via a homomorphism $\sigma:A\to \Aut(X,[\omega])$ is an isomorphism
$$\phi:(V,Q_0)\to (H^k(X,\bbR),Q)$$
such that, for any $g\in A$,
$\rho\circ \rho_0(g)\circ \phi^{-1} = \sigma^*(g)$.
 
Let $f:Y\to T$ be a family of polarized K\"ahler manifolds. Suppose that the group $A$ acts holomorphically on $Y$ preserving the Kahler form $[\omega]$ on $Y$ which defines the polarization form $[\omega]_t$ on each fibre $X_t$. Also we assume that $g(X_t) = X_t$ for all $t\in T$.  Let $V = H^k(X_{t_0},\bbR)$ with polarization form $Q_{t_0}$. Fix an isomorphism $(H^k(X_{t_0},\bbR),Q_{t_0}) \to (V,Q)$ and define $\rho:A\to \Aut(V,Q)$ as the representation $g\mapsto g^*$ of $A$ on $H^k(X_{t_0},\bbR)$. One can show that the monodromy representation $\pi(T,t_0)\to \Aut(V,Q)$ commutes with the representation $\rho$. Let $\Gamma_\rho(f)$ be the centralizer of the group $\rho(A)$ in $\Aut(V,Q)$. Since 
$\Gamma_\rho(f)$ preserves a lattice in  $V$, the image of $H^k(X,\bbZ)$, 
the group $\Gamma_\rho(f)$ is a discrete subgroup of $\Aut(V,Q)$. So the period map
\begin{equation}
\phi:T\to \Gamma_\rho(f)\backslash \calD_\bff(V,Q_0)^\rho.
\end{equation}
is a holomorphic map of complex varieties.   Consider the composition of the multi-valued period map $T\to \calD_\bff(V,Q)^\rho$ with the projection to $\Fl(\bff(\chi),V_\bbC(\chi))^0$. The group $ \Gamma_\rho(f)$ leaves each subspace $V_\bbC(\chi)$ invariant and preserves the flags $(F^p_\chi)$. Since, in general, it is not a discrete subgroup of $\Aut(V_\bbC(\chi),Q)$, only the multi-valued period map is defined as a holomorphic map. By passing to the universal cover $\tilde{T}$ we have a one-valued holomorphic map
\begin{equation}\label{pps}
\phi_\lambda:\tilde{T}\to \Fl(\bff(\chi),V_\bbC(\chi))^0.
\end{equation}
We call this map the {\it eigenperiod map} of the family $f$.

Consider the case when $k = 1$. A Hodge isometry of $V$ is an automorphism of $V$ which is an automorphism of $(V,I)$, where $I$ is a positive complex structure defined by a $Q$-polarized Hodge structure on $V$. Conversely any complex automorphism of $V_\bbC$ preserving the Hodge decomposition and $Q$ arises from a $\bbR$-linear automorphism of $V$ by extension of scalars.  
 
Given a representation $A\to \Sp(V,Q):=\Aut(V,Q)$ we want to describe the set of positive complex structures $I$ on $V$ such that any $g\in A$ commutes with $I$.
In other words, the group $A$ acts naturally on the Siegel half-plane $G(g,V_\bbC)_H$ via its action on $V_\bbC$ and we want to describe the subset 
$$G(g,V_\bbC)_H^\rho = \{W:g(W)\subset W, \forall g\in A\}$$ 
of fixed points of $A$. 

Let $W_\chi = V_\chi \cap W$. Then 
$$V_\bbC(\chi) = W_\chi\oplus \overline{W}_{\chi}.$$
We know that  $H|W_\chi > 0 $ and $H|\overline{W}_{\chi} < 0$. Also  $W$ and $\overline{W}$ are orthogonal with respect to $H$. Therefore the signature of $H|V_\bbC(\chi)$ is equal to  $ (p_\lambda,q_\lambda)$, where $p_\lambda = \dim W_\chi, q = \dim \overline{W}_{\chi}$.

For any  hermitian complex vector space $(E,h)$ of signature $(p,q)$  let us denote by $G(p,(E,h))$ the open subset of the Grassmannian $G(p,E)$ which consists of $p$-dimensional subspaces $L$ such that $h|L > 0$. It is known that $E$ admits a basis $e_1,\ldots,e_{p+q}$  such that the matrix of $E$ is equal to the matrix $I(p,q)$ \eqref{ipq}.
This shows that the group $U(E,h)$ of isometries of $(
E,h)$ is isomorphic to the group
\begin{equation}
U(p,q) = \{A\in \GL(p+q,\bbC): {}^tA\cdot I(p,q)\cdot \bar{A} = I(p,q)\}.
\end{equation}
By Witt's theorem, the group $U(E,h)$ acts transitively on $G(p,(E,h))$ with the isotropy subgroup of $L$ equal to 
$U(L, h|L)\times U(L^\perp, h|L^\perp).$ Thus 
$$G(p,(E,h)) \cong U(p,q)/U(p)\times U(q).$$
It has a structure of a hermitian symmetric space of non-compact type $I$ in Cartan's classification.

Identifying $E$ with $\bbC^g$ and $h$ with the matrix $I(p,q)$, we identify a subspace from $G(p,(E,h))$ with the column space of a matrix $\Pi$ of size $(p+q)\times p$. Writing the matrix in the form
$\Pi = \small{\begin{pmatrix}A\\B\end{pmatrix}},$ where $A$ is a square matrix of size $p$, the positivity condition is expressed by the condition
$$\begin{pmatrix}{}^tA&{}^tB\end{pmatrix}\cdot \begin{pmatrix}I_p&0_{pq}\\
0_{qp}&-I_q\end{pmatrix}\cdot \begin{pmatrix}\bar{A}\\\bar{B}\end{pmatrix} = {}^tA\cdot \bar{A}-{}^tB\cdot \bar{B} > 0.$$
One can show that $|A|\ne 0$ so, replacing $\Pi$ with  $\Pi\cdot A^{-1}$ we may assume that 
$\Pi$ is of the form $\small{\begin{pmatrix}I_p\\Z\end{pmatrix}}$ for some complex $(q\times p)$-matrix. The condition on $Z$ is
\begin{equation}\label{bd2}
I_{p} - {}^tZ\cdot \bar{Z} > 0.
\end{equation}
This shows that $G(p,(E,h))$ is isomorphic to a bounded domain $\calI_{pq}$ in $\bbC^{pq}$ defined by the previous inequality. Any hermitian symmetric space of non-compact type $I$ is isomorphic to such a domain.

Taking  $p=1$ we see that $\calI_{1q}$ is isomorphic as a complex manifold to the unit complex ball of dimension $q$
\begin{equation}
\bbB_q = \{z\in \bbC^q:|z_1|^2+\ldots+|z_q|^2 < 1\}.
\end{equation}
In particulat we we have
$$U(1,q)/U(1)\times U(q) \cong  \bbB_q.$$

There is a natural  embedding of $\calI_{pq}$  into the Siegel half-plane $\calZ_{p+q}$ for which we choose a boundary model from \eqref{bd2}. It is given by assigning to a matrix $Z$ of size $p\times q$ satisfying \eqref{bd1} the symmetric square matrix of size $p+q$ 
\begin{equation}\label{satake}
Z' = \begin{pmatrix}0_p&{}^tZ\\
Z&0_q\end{pmatrix}\end{equation}
satisfying \eqref{bd2}.
Similarly, we have an embedding 
\begin{equation}\label{inclusion}
\prod_{i=1}^k\calI_{p_iq_i}\hookrightarrow \calZ_g,\end{equation}
where $g = p_1+\ldots+p_k = q_1+\ldots+q_k$.
Also we see that the Siegel half-plane $\calZ_g$ is isomorphic to a closed subvariety of $\calI_{gg}$ defined by the equation ${}^tZ = Z$.

Fix  a partition of $\wt(\rho)$.  Let $H_\chi = H|V_\bbC(\chi)$. For any $\chi\in \wt(\rho)^+$, let 
$(p_\chi,q_\chi)$ be the signature of $H_\chi$. By \eqref{real},
the signature of $H_{\bar{\chi}} = (q_\chi,p_\chi)$. The set of pairs of numbers $(p_\chi,q_\chi)$, well-defined up to permutation, is called the {\it type of the representation} $\rho$.

We construct the inverse of the map \eqref{product}
\begin{equation}\label{isom}
\prod_{\chi\in \wt(\rho)^+} G(p_\chi,(V_\bbC(\chi),H_\chi)) \to G(g,V_\bbC)_H^\rho.
\end{equation}
It is defined by assigning to  a collection of subspaces $(E_\chi)$ from the left-hand-side the direct sum
$$W = \bigoplus_{\chi\in \wt(\rho)^+} E_\chi \oplus\overline{(E_\chi)_{H_\chi}^\perp}.$$
It is clear that $\overline{(E_\chi)_{H_\chi}^\perp} \subset V_{\bar{\chi}}$ and the restriction of $H$ to this subspace is positive. Thus $H|W > 0$ and 
$$\dim W = \sum_{\chi\in \wt(\rho)^+} p_\chi+ \sum_{\chi\in \wt(\rho)^-} q_\chi = g.$$
 Also, using \eqref{orth} and the fact that 
$H(E_\chi,(E_\chi)_{H_\chi}^\perp) = 0$ implies $Q(E_\chi,\overline{(E_\chi)_{H_\chi}^\perp}) = 0$, we see that $Q|W = 0$. Thus $W\in G(g,V_\bbC)_H$. As a sum of eigensubspaces it is obviously $\rho$-invariant, so the map is well-defined and is easy to see coincides with the inverse of the map \eqref{product}.

Thus for any linear  representation $\rho:A\to \Sp(V,Q)$ with 
$\wt(V_\bbC) = \{\chi_i, \bar{\chi}_i, i = 1,\ldots,k\}$  and $(p_i,q_i) = (p_{\chi_i},q_{\chi_i})$, we obtain an embedding of hermitian symmetric domains
$$\prod_{i=1}^k\calI_{p_i,q_i} \hookrightarrow \calZ_g.$$

\begin{ex}
Let $A  =(g)$ be a cyclic group of order $4$. Assume that its generator $g$ satisfies $g^2 = -1$. Therefore its action on $V$ is equivalent to a complex structure  $I$ on $V$ and hence defines a subspace $L\subset V_\bbC$ on which $g$ acts with eigenvalue $i$. Thus 
$$L = (V_\bbC)_\chi,\quad  \bar{L} = (V_\bbC)_{\bar{\chi}},$$
where $\chi(g) = i, \bar{\chi}(g) = -i.$
Let $Q$ be a symplectic form on $V$ such that $Q(v,I(v'))$ is a symmetric form of signature $(p,q)$ and let $H= iQ(x,\bar{x})$ be the associated hermitian form. Its restriction $h$ to $L$ is of signature $(p,q)$. Let 
$E$ be a positive subspace of $L$ of dimension $p$. Let $E_h^\perp$ be its orthogonal complement in $L$ with respect to $h$. The restriction of $h$ to $E_h^\perp$ is negative. Its conjugate subspace $E' = \overline{E_h^\perp}$ belongs to $\overline{W}$. Since
$H(x,\bar{x}) = -H(\bar{x},x)$,  the restriction of $H$ to $E'$ is positive. Thus  $W = E\oplus E'$ defines a positive complex structure on $V$ with respect to $Q$. The operator $I$ acts on a $p$-dimensional subspace $E$ of $W$ as $i1_E$ and acts on the orthogonal $q$-dimensional subspace $E'$ as $-i1_{E'}$. The decomposition $V_\bbC = W\oplus \overline{W}$ is a $Q$-polarized Hodge structure invariant with respect to the representation of $A$ in $V$.

In coordinates, choose a standard symplectic basis $e_i$ in $V$ for a symplectic form $Q$, and define $I$ by the formula 
$$I(e_k) = e_{k+g}, \ \ k = 1,\ldots,p,\quad   I(e_{p+k}) = -e_{p+g+k}, \ k = 1,\ldots,q.$$ Then
$Q(v,I(v'))$ is a symmetric quadratic form of signature $(p,q)$. 
Then $I$ preserves the lattice in $V$ spanned by the vectors $e_i$ and acts on $W$ preserving the projection $\Lambda$ of this lattice in $W$.  Thus $g$ acts an automorphism of 
the principally polararized abelian variety $A = W/\Lambda$. This automorphism is of type $(p,q)$, i.e.  $W$ decomposes into a direct sum of eigensubspaces of dimension $p$ and $q$.

We can find a basis in $V_\bbC$ which is a sum of a basis in $L$ and $\bar{L}$ such that $W = E\oplus E'$ is represented by a matrix
$$\begin{pmatrix}I_p&0_{pq}\\
Z&0_q\\
0_q&I_{q}\\
0_{qp}&{}^tZ\end{pmatrix}$$
where $Z\in \calI_{pq}$. After an obvious change of a basis this matrix becomes a matrix \eqref{satake} defining the embedding of $\calI_{pq}$ in $\calZ_g$.

This shows that the isomorphism classes of principally polarized abelian varieties of dimension $g$ admitting an automorphism of order 4 of type $(p,q)$ is isomorphic to the space of orbits
$$\Gamma_{pq}\backslash \calI_{pq} \subset \calA_g = \Sp(2g,\bbZ)\backslash\calZ_g ,$$ where $\Gamma_{pq}$ is the subgroup of matrices $M\in \Sp(2g,\bbZ)$ which commute with the matrix 
$\begin{pmatrix}0_g&I_{p,q}\\
-I_{p,q}&0_g\end{pmatrix}$ defining the operator $I$. 

\end{ex}

\section{Arrangements of hyperplanes}\label{sec7} 
Let $H_1,\ldots,H_{m}$ be hyperplanes  in $\bbP^n(\bbC)$ defined by linear forms 
$$f_i(t_0,\ldots,t_n) = \sum_{i=0}^na_{ij}t_j, \quad i = 1,\ldots,m.$$
 We assume that the hyperplanes are in a general position. This means that all maximal minors of the matrix
\begin{equation}\label{mat}
M= (a_{ij})\end{equation} are not equal to zero. Geometrically this means that the intersection of any $n+1$ hyperplanes is a point. Let  
 $$U =  \bbP^n\setminus \bigcup_{i=1}^m H_i$$ be the complementary set. Choose a set of rational numbers 
$$\boldsymbol{\mu} = (\mu_1,\ldots,\mu_m),$$
satisfying
\begin{equation}\label{cond1}
0 < \mu_i < 1;\end{equation}
\begin{equation}\label{cond2}
|\boldsymbol{\mu}|: =\sum_{i=1}^m \mu_i\in \bbZ.\end{equation}

It is well-known that the fundamental group $\pi_1(U,u_0)$ is a free group generated by the homotopy classes $g_i$ of loops which are defined as follows. Choose a line $\ell$ in $\bbP^n$ which does not intersect any codimension 2 subspace $H_i\cap H_j$ and contains $u_0$.  Let $\gamma_i$ be a path in $\ell$ which connects $u_0$ with  a point on a small circle around the point $p_i = \ell\cap H_i$, goes around the circle in a counterclockwise way, making the full circle, and then returns in the opposite direction to the point $u_0$. The homotopy classes $g_i =[\gamma_i]$ generate the group $\pi_1(U,u_0)$ with the defining equation $g_1\cdots g_{m} = 1$. 

Let $\calL_{\boldsymbol{\mu}}$ be the complex local coefficient system on $U$ defined by the homomorphism
\begin{equation}\label{ch}
\chi:\pi_1(U,u_0)\to \bbC^*, \ g_i\mapsto e^{-2\pi i\mu_i}.
\end{equation}

 Let $d$ be the least common denominator of the numbers $\mu_i$ and 
 $$A_d(m) = (\bbZ/d)^{m}/\Delta(\bbZ/d),$$
where $\Delta$ denotes the diagonal map.   Consider the finite cover $U'$ of $U$ corresponding to the homomorphism 
$\pi_1(U,u_0)\to A_d(m) $ defined by sending $g_1$ to the vector $e_1 = (1,0,\ldots,0)$, and so on. Let $X$ be the normalization of $\bbP^n$ in the field of rational functions on $U'$. This is a nonsingular algebraic variety  which can be explicitly described as follows.  Consider the map 
$$f:\bbP^n\to \bbP^{m-1}, \quad (t_0,t_1,\ldots,t_n) \mapsto (f_1(t),\ldots,f_m(t)).$$ 
Its image is a linear $n$-dimensional subspace $L$ of $\bbP^{m-1}$. Let $\phi:\bbP^{m-1}\to \bbP^{m-1}$ be the cover defined by 
$$(z_0,\ldots,z_{m-1}) \to (z_0^d,\ldots,z_{m-1}^d).$$
Then $X$ is isomorphic to the pre-image of $L$ under the cover $\phi$.  If we write the linear forms 
$f_i, i > n+1, $ as linear combinations $\alpha_{i0}f_{1}+\ldots +\alpha_{in}f_{n+1}$, then $X$ becomes isomorphic to the subvariety  in $\bbP^{m-1}$ given by the equations
$$-x_{i}^d+\alpha_{i0}x_0^d+\ldots+\alpha_{in}x_n^d = 0,\quad i = n+1,\ldots,m.$$
The group $A_d(m)$ acts naturally on $X$ via multiplication of the coordinates $x_i$ by $d$th roots of unity. The quotient space $X/A$ is isomorphic to $\bbP^n$. The group $A_d(m)$ acts on $V = H^n(X,\bbR)$ and on its compexification $V_\bbC = H^n(X,\bbC)$. 
 
\begin{lem}\label{lem5} 
Let $\chi\in \hat{A}_d(m)$ be defined by \eqref{ch}. Then 
$$H^j(U,\calL_{\boldsymbol{\mu}}) = 0, j\ne n,\quad H^n(U,\calL_{\boldsymbol{\mu}})  \cong H^n(X,\bbC)(\chi).$$
\end{lem}

\begin{proof} Let
$$U' = X\setminus \bigcup_{i=0}^{m-1} \{x_i = 0\}$$
and let $\pi:U'\to U$ be the natural projection. This is a unramified Galois covering with the Galois group $A_d(m)$. 
The direct image of the constant coefficient system $\bbC_{U'}$ in $U$ decomposes into the direct sum of local 
coefficients systems $\calL_\chi$ corresponding to different characters $\chi\in \hat{A}_d(m)$
$$\pi_*(\bbC_{U'}) = \bigoplus_{\chi\in \hat{A}} \calL_\chi.$$
We have an isomorphism of cohomology with compact support
$$H_c^j(U,\pi_*(\bbC_{U'})) \cong H_c^j(U',\bbC)$$
which is compatible with the action of $A$ and gives an isomorphism
$$H_c^j(U',\bbC)(\chi) \cong H_c^j(U,\calL_{\boldsymbol{\mu}}).$$
One can show (see \cite{DM}, p. 20) that for any character $\chi$ determined by the numbers $\mu_i$ satisfying \eqref{cond1} and \eqref{cond2}, 
$H_c^j(U,\calL_\chi) \cong H^j(U,\calL_\chi)$. Let $Y = X\setminus U'$. Now consider the long exact sequence
$$\ldots \to H_c^{j}(U',\bbC) \to H^j(X,\bbC) \to H^j(Y,\bbC) \to H_c^{j+1}(U,\bbC)\to \ldots.$$
Since $U'$ is affine, it has homotopy type of a CW-complex of dimension $n$, and hence $H_c^j(U',\bbC) = 0$ for $j > n$.
$H^{n-1}(Y,\bbC)$ is generated by the fundamental classes of the irreducible components 
$Y_i = Y\cap \{x_i = 0\}$. They are fixed under the action of $A_d(m)$, and hence  $H^{n-1}(Y,\bbC)_\chi = 0.$ Since $Y_i$ and $X$ are complete intersections in projective space, by Lefschetz's Theorem on a hyperplane section, the natural homomorphism
$H^{j}(\bbP^{m-1},\bbC) \to H^j(X,\bbC),$ is an isomorphism for $j < n$. Similarly,
$H^{j}(\bbP^{m-2},\bbC) \to H^j(Y_i,\bbC),$ is an isomorphism for $j < n-1$. Since the group $A_d(m)$ acts trivially on cohomology of the projective space, we obtain that 
$$ H^j(X,\bbC)(\chi) = \{0\},\quad  j \ne n .$$
Also using the Mayer-Vietoris sequence, one shows that 
$$  H^j(Y,\bbC)(\chi) = 0, \quad  j \ne n-1.$$
 The long exact sequence now gives
$$H^j(U,\calL_{\bldmu}) \cong H_c^j(U',\bbC)(\chi) = 0, \quad j \ne n, $$
$$ H^n(U,\calL_{\bldmu}) \cong H_c^n(U',\bbC)(\chi) \cong  H^n(X,\bbC)(\chi). $$
\end{proof}

\begin{lem}\label{my} 
For any character $\chi$ defined by a collection $\boldsymbol{\mu}$,
$$\dim  H_\chi^{pq}(X) = \binom{|\boldsymbol{\mu}|-1}{p}\binom{m-1-|\boldsymbol{\mu}|}{q},$$
$$\dim H^n(X,\bbC)_\chi = \binom{m-2}{n}.$$
\end{lem}

\begin{proof} We only sketch the proof. Another proof was given by P. Deligne \cite{De}. One can compute explicitly the Hodge decomposition of a nonsingular complete intersection $n$-dimensional subvariety $X$ in projective space $\bbP^N$.  Let 
$F_1 = \ldots = F_{N-n} = 0$ be homogeneous equations of $X$ in variables $x_0,\ldots,x_N$. Let $y_1,\ldots,y_{N-n}$ be new variables and
$$F = y_1F_1+\ldots+y_{N-n}F_{N-n}.$$
Define the {\it jacobian algebra} of $F$ as the quotient algebra
$$R(F) = \bbC[x_0,\ldots,x_N,y_1,\ldots,y_{N-n}]/J_F,$$
where $J_F$ is the ideal generated by partial derivatives of $F$ in each variable. It has a natural bigrading defined 
by $\deg x_i = (1,0), \deg y_j = (0,1)$. There is a natural isomorphism (see \cite{Ter})
$$H_{prim}^{pq}(X) \cong R(F)_{q,((N+1+q)d-n-1)},$$
where $d$ is the sum of the degrees of the polynomials $F_i$. If $G$ is a group of automorphisms of $X$ induced by linear transformations of $\bbP^N$, then its action on the cohomology is compatible with the action on the ring $R(F)$ (the representation on this ring must be twisted by the one-dimensional determinant representation). In our case the equations $f_i$ are very simple, and the action of the  group $A_d(m)$  can be explicitly computed. 

The second formula follows from the first by using the combinatorics of binomial coefficients.
\end{proof} 

Let $\chi$ be the character corresponding to $\boldsymbol{\mu}$. Using the isomorphism from Lemma \ref{lem5}, we can define
\begin{equation}
H^0(U,\calL_{\boldsymbol{\mu}})^{pq} = H_\chi^{pq}(X).\end{equation}
One has the following description of this space (see \cite{DM},\cite{De}).
$$H^0(U,\calL_{\boldsymbol{\mu}})^{pq} = H^q(\bbP^n,\Omega_{\bbP^n}^p(\log \cup H_i)\otimes \calL_{\boldsymbol{\mu}}),$$
where $\Omega_{\bbP^n}^p(\log \cup H_i)$ is the sheaf of meromorphic differential $p$-forms with at most simple poles on the hyperplanes $H_i$. These forms can be considered as multi-valued forms on $\bbP^n$ with local branches defined by the local coefficient system. 

Consider the family $\calX$ parametrized by the set of all possible ordered sets of $m$ hyperplanes in general linear position. Any such collection is defined uniquely by the matrix $M$ \eqref{mat} of coefficients of linear forms defining the hyperplanes. Two matrices define the same collection if and only if their columns are proportional. The set of equivalence classes is an algebraic variety $\calX_{n,m}$ of dimension $mn$. Let $\tilde{\calX}_{n,m}$ be the universal covering of $\calX_{n,m}$. We have the eigenspace period map
$$\phi:\tilde{\calX}_{n,m} \to \Fl(\bff(\chi),V_\bbC(\chi))^0,$$
where $V = H^n(X_{t_0},\bbR)$ for some $t_0\in \calX_{n,n}$ and 
$$f_p(\chi) = \sum_{p'\ge p} \binom{|\boldsymbol{\mu}|-1}{p'}\binom{m-1-|\boldsymbol{\mu}|}{n-p'}, \ p = 0,\ldots, n.$$

The group $\GL(n+1,\bbC)$ acts naturally on $\calX_{n,m}$ by left-multiplications. Geometrically this means a projective transformation sending a collection of hyperlanes to a projectively equivalent collection. On can show that the orbit space $\calP_{n,m} = \calX_{n,m}/\GL(n+1,\bbC)$ is a quasi-projective algebraic variety of dimension $(m-n-1)m$. Two $X_t$ and $X_{t'}$ corresponding to points in the same orbit are $A_d(m)$-equivariantly isomorphic. This shows that the eigenperiod map is defined on the universal covering $\tilde{\calP}_{n,m}.$

The following result is a theorem of A. Varchenko \cite{Va}.

\begin{thm}\label{var} Assume that $\chi$ corresponds to collection of numbers $\mu_i$ such that
$$|\boldsymbol{\mu}| = \mu_1+\ldots+\mu_m = n+1,$$
or, equivalently, $h_\chi^{n0}(X) = 1$.  Then the eigenperiod map
$$ \tilde{\calP}_{n,m}\to  \Fl(\bff(\chi),V_\bbC(\chi))^0$$
is a local isomorphism onto its image.
\end{thm}

\begin{ex} We take $n = 1$. In this case $U = \bbP^1-\{p_1,\ldots,p_m\}$. Without loss of generality we may assume that $p_{m}$ is the infinity point and other $p_i$'s have affine coordinates $z_i$. 

 Consider a multi-valued form
$$\omega =g(z)(z-z_1)^{-\mu_1}\ldots (z-z_{m-1})^{-\mu_{m-1}}dz,$$
where $g(z)$ is a polynomial of degree $|\boldsymbol{\mu}|-2$ (zero if $|\boldsymbol{\mu}| < 2$). Under the monodromy transformation the power $(z-z_i)^{-\mu_i}$ is transformed in the way prescribed by the local coefficient system. At the infinity we have to make the variable change $z = 1/u$ to transform $\omega$ to the form
$$\omega =h(u)(1-z_1u)^{-\mu_1}\ldots (1-z_{m-1}u)^{-\mu_{m-1}}u^{-\mu_m}du.$$
Again we see that the local mondromy $u\mapsto  ue^{2\pi i\phi}$  agrees with the mondromy of the local coefficient system. This shows that the forms  $\omega$ span a subspace  of  dimension $|\boldsymbol{\mu}|-1$. Comparing with the formula from Lemma \ref{my} we see that they span the whole space $H^1(U,\calL_{\boldsymbol{\mu}})^{10} = H_\chi^{10}(X).$

Assume 
\begin{equation}\label{2}
|\boldsymbol{\mu}|  = 2,\end{equation}
or, equivalently, $h_\chi^{10}(X) = 1$. In this case, the space $F_\chi^1 = H^1(U,\calL_{\bldmu})^{10}$ is generated by the form
$$\omega_{\boldsymbol{\mu}} =(z-z_1)^{-\mu_1}\ldots (z-z_{m-1})^{-\mu_{m-1}}dz.$$
The form $\omega_{\boldsymbol{\mu}}$ is called the hypergeometric form with exponents $(\mu_1,\ldots,\mu_{m-1})$. 
The space $V = H^1(X,\bbR)(\chi)$ is of dimension $m-2$ and $\Fl(\bff(\chi),V_\bbC(\chi))^0$ is the open subspace of the projective space of lines in $V_\bbC$ such that restriction of the hermitian form 
$iQ(x,\bar{x})$ to $F^1$ is positive definite. Thus $\Fl(\bff(\chi),V_\bbC(\chi))^0$ is isomorphic to a complex ball $\bbB_{m-3}$ and we have  the eigenperiod map
\begin{equation}\label{ep1}
\tilde{\calP}_{1,m} \to \bbB_{m-3}.
\end{equation}
 Note that 
$\dim \calP_{1,m} = m-3 = \dim \bbB_{m-3}$. By Theorem  \ref{var}, this map is a local isomorphism (in the case $n = 1$ this is an earlier result which can be found in \cite{DM}.)

Now recall some constructions from the Geometric Invariant Theory. Consider the space $X = (\bbP^n)^m$ parametrizing $m$-tuples of points in $\bbP^n$, not necessary distinct. Consider a map 
$f:(\bbP^n)^m \to \bbP^N$ given by all multi-homogenous monomials  of multi- degrees 
$\bfk = (k_1,\ldots,k_m), \ k_i > 0$. The group $\SL(n+1)$ acts naturally on $\bbP^n$ by projective transformations and on the product by the diagonal action. Let $S$ be the projective coordinate ring of the image of the product in $\bbP^N$. The group $G = SL(n+1,\bbC)$ acts on $S$ via its linear action on $\bbP^N$. Let $S^G$ be the subring of invariant elements. By a theorem of Hilbert, the subring $S^G$ is finitely generated graded algebra over $\bbC$. By composing the map $f$ with the Veronese map $\bbP^N\to \bbP^{M}$ defined by all homogeneous monomials of some sufficiently large degree one may assume that $S^G$ is generated by elements of degree 1. A standard construction realizes this ring as the projective coordinate ring of some projective variety. This projective variety (which is uniquely defined up to isomorphism) is denoted by 
$X/\!/_\bfk G$ and is called the $GIT$-quotient of $X$ by $G$ with respect to the linearization defined by numbers $\bfk$. Assume $n = 1$. Let $U^{\s}$ (resp. $U^{\ses}$) be the open Zariski subset of $(\bbP^1)^m$ parametrizing collections of points  $(p_1,\ldots,p_m)$ such that for any $p\in \bbP^1$ the following conditions is satisfied:
\begin{equation}\label{stab}
2\sum_{p_i =p} k_i < \sum_{i=1}^m k_i \ (\text{resp. $2\sum_{p_i =p} k_i \le \sum_{i=1}^m k_i$)}.
\end{equation}
One shows that there is a canonical surjective map $U^{\ses} \to X/\!/_\bfk G$ whose restriction to $U^{\s}$ is isomorphic to the natural projection of the orbit space $U\to U/G$. Two point sets in the set $U^{\ses}\setminus U^s$ are mapped to the same point if and only if they have the same subset of points $p_j, j\in J$, satisfying $2\sum_{j\in J}k_i =  \sum_{i=1}^m k_i$ or its complementary set.

Let $P_{1,\bfk}^{\ses}$ (resp. $P_{1,\bfk}^{\s}$) denote the GIT-quotient $(\bbP^1)^m/\!/_\bfk \SL(2,\bbC)$
(resp. its open subset isomorphic to the orbit space  $U^{\s}/\SL(2)$). Our variety  $P_{1,m}$ is isomorphic to an open subset of $P_{1,\bfk}^{\s}$ for any $\bfk = (k_1,\ldots,k_m)$.

The following is a main result of Deligne-Mostow's paper \cite{DM}.

\begin{thm}\label{dm} Let $\bldmu$ satisfies condition \eqref{2}. Write $\mu_i = k_i/d$, where $d$ is a common least denominator of the $\mu_i$'s. Assume that the following condition is satisfied:
\begin{equation}\label{int2} 
(1-\mu_i-\mu_j)^{-1} \in \bbZ,  \quad \text{for any $i\ne j$ such that $\mu_i+\mu_j < 1$}.\end{equation}
Then the monodromy group $\Gamma_{\bldmu}$ is a discrete group of holomorphic  automorphisms of $\bbB_{m-3}$ and the eigenperiod map \eqref{ep1} extends to an isomorphism
$$\Phi:P_{1,\bfk}^{\s} \to \Gamma_{\bldmu}\backslash \bbB_{m-3}.$$
Moreover, $\Phi$ can be extended a holomorphic isomorphism 
$$\overline{\Phi}: P_{1,\bfk}^{\s} \to \overline{\Gamma_{\bldmu}\backslash \bbB_{m-3}},$$
where the target space is a certain compactification of the ball quotient obtained by adding a finite set of points. 
\end{thm}

The list of collections of numbers $(\mu_1,\ldots,\mu_m)$ satisfying \eqref{int2} is finite for $m > 4$.  It consists of 27 collections with $m = 5$, of 7 collections with $m = 6$ and one collection with $m = 7$ and $m = 8$. One can extend the theorem to the case when $\bldmu$ satisfies a weaker condition:
\begin{equation}\label{mostow}
(1-\mu_i-\mu_j)^{-1}\in \bbZ \quad \text{for any $i\ne j$ such that $\mu_i+\mu_j < 1, \mu_i\ne \mu_j$}
\end{equation}
$$2(1-\mu_i-\mu_j)^{-1}\in \bbZ \quad \text{for any $i\ne j$ such that $\mu_i+\mu_j < 1, \mu_i = \mu_j$}.$$
This gives additional cases, with largest $m$ equal to 12.

 Write $\mu_i = k_i/d$ as in Theorem \ref{dm}. Recall 
that, by Lemma \ref{lem5} we have an isomorphism
$$H^1(U,\calL_{\boldsymbol{\mu}}) \cong H^1(X,\bbC)(\chi),$$
where $X$ is the curve in $\bbP^{m-1}$ defined by  equations
$$-x_i^d+a_{i0}x_0^d+a_{i1}x_1^d = 0,\ i = 2,\ldots,m-1.$$
Let $H = \Ker(\chi)$ and $X_{\boldsymbol{\mu}} = X/H$. one can show that the curve $X_{\bldmu}$ is isomorphic to a nonsingular model of the affine algebraic curve with equation
\begin{equation}\label{eq}
y^d =  (x-z_1)^{k_1}\ldots (x-z_{m-1})^{k_{m-1}}.\end{equation}
The cyclic group $C_d = \bbZ/d\bbZ$ acts naturally on this curve by acting on the variable $y$. One proves that 
$$H^1(U,\calL_{\bldmu}) \cong H^1(X_{\bldmu},\bbC)(\chi).$$
For example, take $m = 6, d = 3, \mu_1= \ldots=\mu_6 = 1$. The curve $X_{\bldmu}$ is the curve of genus 4 curve with affine equation $y^3 = (x-z_1)\cdots(x-z_6).$ We have $\dim H^{10}(X_{\bldmu})(\chi) = 1$. Thus the action of $C_3$ on the Jacobian variety is of type $(1,3)$. As we saw in section 
\ref{ab}, the locus of principally polarized abelian 4-folds admitting an automorphism of order 3 of type $(1,3)$ is isomorphic to a 3-dimensional ball quotient. This agrees with the theory of Deligne-Mostow. 

\end{ex}

\section{K3 surfaces}

This is our second example. This time we consider  compact orientable simply-connected 4-manifolds $M$. Recall that by a theorem of M. Friedman, the homeomorphism type of $M$ is determined by the cup-product on $H_2(M,\bbZ) \cong \bbZ^{b_2(M)}$. It is a symmetric bilinear form whose matrix in any $\bbZ$-basis has determinant $\pm 1$. The corresponding integral quadratic form is called in such case a {\it unimodular quadratic form}.  The group $H_2(M,\bbZ)$ equipped with the cup-product is an example of a {\it lattice}, a free abelian group equipped with a non-degenerate  symmetric bilinear form, or, equivalently, with an $\bbZ$-valued quadratic form.  We will denote the values of the bilinear form by $(x,y)$ and use $x^2$ to denote $(x,x)$. We apply the terminology of real quadratic forms to a lattice whose quadratic form is obtained by extension of scalars to $\bbR$. There are two types of indefinite lattices. The {\it even type} is when the quadratic form takes only even values, and the {\it odd type} when it  takes some odd values. A unimodular indefinite lattice  of odd type can be given in some $\bbZ$-basis  by a matrix $I(p,q)$ \eqref{ipq}.
Its Sylvester signature is $(p,q)$. A unimodular indefinite quadratic form $S$ of even type of signature $\le 0$ is isomorphic to the orthogonal sum $U^{\oplus p}\oplus E_8^{\oplus q}$, where $U$ is given by the matrix 
$$
\begin{pmatrix}0&1\\
1&0\end{pmatrix}$$
and $E_8$ is a unique (up to isomorphism) even negative definite unimodular lattice of rank 8. 
The index of $S$ is equal to $-8q$ and the Sylvester signature is equal to $(p,p+8q)$. The lattice $E_8$ is defined by the Dynkin diagram of type $E_8$. Its vertices correspond to vectors with $(x,x) = -2$ and $(x_i,x_j) = 1$ or $0$ dependent on whether the corresponding vertices are joined by an edge. Note that changing the orientation changes the index $I(M)$ to its negative. Also by a theorem of Donaldson $p\ne 0$ for the even cup-product on a simply-connected 4-manifold.

A 4-manifold $M$ is called a {\it $K3$-manifold} if it realizes the form $L_{K3} = U^3\oplus E_8^2$, where the direct sum is considered to be an orthogonal sum.

\begin{thm} Each $K3$-manifold is homeomorphic  to a compact simply connected complex K\"ahler surface $X$ 
with trivial first Chern class $($i.e. the second exterior power of the holomorphic tangent bundle is the trivial line bundle$)$. 
\end{thm}

\begin{proof} It is enough to construct such a complex surface realizing the quadratic form $L_{K3}$. 
Let $X$ be a nonsingular quartic surface in $\bbP^3(\bbC)$. The standard exact sequence
$$0\to T_X \to T_{\bbP^3}|X\to N_{X/\bbP^3}\to 0$$
allows one to compute the Chern classes of $X$. We have 
$$1+c_1+c_2 = (1+h)^4/(1+4h) = 1+6h^2,$$
where $h$ is the class of a plane section. Thus $c_1 = 0$ and $c_2 = 6\cdot 4 = 24.$ Since $c_1 = 0$, the canonical class $K_X = -c_1 = 0$. Thus $H^{20}(X) =H^0(X,\Omega_X^2) = \bbC$. Since $24 = c_2 = \chi(X) = \sum (-1)^ib_i(X)$ and $X$ is simply-connected, we get $b_2(X) = 22$.  Thus $b_2 = h^{20}+h^{11}+h^{02}$ implies 
$h^{11} = 20$. By the Index formula \eqref{index} we obtain that $I(M) = -16$, hence the Sylvester signature is $(3,19)$. One also uses Wu's Theorem which asserts that the cup-product is even if $c_1(X)$ is divisible by  2 in $H^2(X,\bbZ)$. Thus $H^2(X,\bbZ)$ is an even lattice. This implies that the cup-product $H^2(X,\bbZ)$ is isomorphic to the lattice $L_{K3}.$

\end{proof}

A non-trivial theorem which is a corollary of the theory of periods of K3 surfaces is the following.

\begin{thm} All complex $K3$ surfaces are diffeomorphic to a nonsingular quartic surface in $\bbP^3(\bbC)$.
\end{thm}

Let $X$ be an algebraic K3 surface (e.g. a quartic surface). We have already observed that  $(h^{20},h^{11},h^{02}) = (1,20,1)$. The Hodge flag is 
$$0 \subset F^2 = H^{20}(X) \subset F^1 = H^{20}(X)+H^{11}(X) \subset F^0 = H^2(X,\bbC).$$

The Hodge structure on cohomology $V = H_{prim}^2(X,\bbR)$ is a AHS of weight 2 of type $(1,19,1)$,  where we take the polarization defined by a choice of a Hodge form $h\in H^2(X,\bbZ)$. The polarization form $Q$ admits an integral structure with respect to the lattice $H^2(X,\bbZ)$. Note that the flag $(F^p)$ is completely determined by $F^2$ (because  $F^1 = (F^2)^\perp$). This shows that the period space  $D_\bff(V,Q)$ is isomorphic, as a complex manifold, to 
\begin{equation}
\calD_h(V) = \{\bbC\cdot v\in \bbP(V_\bbC):Q(v,v) = 0, Q(v,\bar{v}) > 0\}.
\end{equation}
 A {\it marking} of a K3 surface is an isomorphism of lattices   $\phi:L_{K3} \to H^2(X,\bbZ)$. After tensoring with $\bbR$ or $\bbC$, it defines an isomorphism 
 $$\phi_\bbR:(L_{K3})_\bbR \to H^2(X,\bbR),\  \phi_\bbC:(L_{K3})_\bbC\to H^2(X,\bbC).$$ Choose a marking. Let $l = \phi^{-1}(h)$. Let $l^\perp$ denote the orthogonal complement of $\bbZ l$ in $L_{K3}$. It is a sublattice of signature $(2,19)$. The period space $\calD_h(V)$ becomes isomorphic to 
 the 19-dimensional complex variety
 $$\calD_{l} = \{\bbC v\in \bbP((l^\perp)_\bbC): (v,v) = 0, (v,\bar{v}) > 0\}.$$
 
 Note that, in coordinates, the period space $\calD_l$ is isomorphic to a subset of lines in $\bbP^{20}$ whose complex coordinates 
 $(t_0,\ldots,t_{19})$ satisfy
 $$t_0^2+t_1^2-t_3^2-\ldots-t_{19}^2 = 0,$$
 $$|t_0|^2+|t_1|^2-|t_3|^2-\ldots-|t_{19}|^2 > 0.$$

One can give another model of the period space as follows. For any line $\bbC v\in \calD_l$ write $v= x+iy$, where $x,y\in (l^\perp)_\bbR$ and consider  the  plane $P$ in $V = (l^\perp)_\bbR$  spanned by $x,y$. We have $(x+iy,x-iy) = x^2+y^2 > 0$. Also  $P$ carries a canonical orientation defined by the orientation class of the basis formed by $x$ and $y$, in this order. This defines a map 
$$\calD_l \to G(2,V)_Q^+ = \{P\in G(2,V)^0: Q|P > 0\}, $$
where $G(2,V)^0$ is the Grassmannian of oriented planes in $V$. It is easy to construct the inverse map. For any $P\in G(2,V)_Q^+$ with a basis $(x,y)$ we assign a line in $V_\bbC$ spanned by $x+iy$. A different basis in the same orientation class changes the complex vector by a complex multiple. This gives a well-defined map to the projective space $\bbP(V_\bbC)$ and, as is easy to see, the image is equal to $\calD_l$.  The changing of orientation in the plane $P$ corresponds to  replacing   $x+iy$ with $x-iy$, i.e. replacing $v$ with $\bar{v}$. If $v$ arises as the period of a K3 surface, this  switching is achieved by the changing the complex structure of the surface to the conjugate one.

The real Grassmannian model shows that the group $\O(V) \cong O(2,19)$ acts transitively on $\calD_l$ with the isotropy group of some point isomorphic to the subgroup $\SO(2)\times \O(19)$. Here $\O(p,q)$ denotes the orthogonal group of the  real quadratic form on $\bbR^n$ defined by the diagonal matrix $I(p,q)$, and $\SO(p,q)$ denotes its subgroup of matrices with determinant 1. So, we obtain
$$\calD_l \cong \O(2,19)/\SO(2)\times \O(19).$$
This description shows that the period space is not connected. It consists of two disconnected copies of a Hermitian symmetric space of non-compact type IV isomorphic to the orbit space
$ \SO(2,19)/\SO(2)\times \SO(19)$.  
In coordinate description of $\calD_l$ the connected components are distinguished by the sign of $\text{Im}(\frac{t_1}{t_0}).$ 

It is known that the K3-lattice $L_{K3}$ represents any even number $2d$, i.e. the set of  elements $x\in L_{K3}$ with $x^2 = 2d$ is not empty. This is easy to see, for example, considering the sublattice of $L_{K3}$ isomorphic to $U$. Less trivial is the fact that any two vectors $x$ and $y$ with $x^2 = y^2$ differ by an isometry of $L_{K3}$ provided that none of them is divisible by an integer $\ne \pm 1$ (i.e. $x,y$ are {\it primitive elements} of the lattice). Let $(X,h)$ be a polarized algebraic K3 surface of degree $h^2 = 2d$.  As always in this case we assume that the polarization class $h$ belongs to $H^2(X,\bbZ)$, i.e. is a Hodge class. We can always choose $h$ to be a primitive element in $H^2(X,\bbZ)$. For each positive integer $d$ fix a primitive element $l\in L_{K3}$ with $l^2 = 2d$. Let $\phi:L_{K3}\to H^2(X,\bbZ)$ be a marking of $(X,h)$. Then $\sigma(\phi^{-1}(h)) = l$ for some $\sigma\in \O(L_{K3})$. This shows that $(X,h)$ always admits a marking such that $\phi^{-1}(h) = l$. We call it a {\it marking of a polarized $K3$ surface}. Its {\it period}
$$p_{X,\phi} = \phi_\bbC^{-1}(H^{20}(X))$$
belongs to  $\calD_l.$ Note that the polarization $h$ is determined by the marking since, by definition, $\phi(l) = h$, so there is no need to use the notation $p_{X,h,\phi}$. 

We say that two marked polarized K3 surfaces $(X,\phi)$ and $(X',\phi')$ are isomorphic if there exists an isomorphism of algebraic varieties $f:X'\to X$ such that $\phi' =  f^*\circ \phi$ and $f^*(h') = h$. Obviously this implies that $p_{X,\phi} = p_{X',\phi'}$. Conversely, suppose this happens. Then 
$\psi=\phi\circ \phi{}^{-1}:H^2(X',\bbZ) \to H^2(X,\bbZ)$ is an isomorphism of lattices such that 
$\psi(H^{20}(X')) = H^{20}(X)$ and $\psi(h') = h$. 

The following is a fundamental result of   I. R. Shafarevich and I.I. Piyatetsky-Shapiro \cite{SPS}, called the {\it Global Torelli Theorem for polarized algebraic $K3$ surfaces}.

\begin{thm} Let $(X,h)$ and $(X',h')$ be two polarized algebraic $K3$ surfaces. Suppose there is an isometry of lattices $\phi:H^2(X',\bbZ) \to H^2(X,\bbZ)$ such that $\phi(h') = h$ and $\phi_\bbC^*(H^{20}(X')) = H^{20}(X)$. Then there exists a unique isomorphism of algebraic varieties $f:X\to X'$ such that $f^* = \phi$. 
\end{thm}

Let 
$$\Gamma_{l} = \{\sigma\in \O(L_{K3}):\sigma(l) = l\}.$$
For any marked polarized K3 surface $(X,\phi)$ with $p_{X,\phi}\in \calD_l$ and any $\sigma\in \Gamma_l$ we have $p_{X,\phi\circ\sigma}\in \calD_l$. So the orbit $\Gamma_l\cdot p_{X,\phi}$ depends only on isomorphism class of $(X,h)$.  The group $\Gamma_l$ is isomorphic to a discrete subgroup of $\O(2,19)$ which acts transitively on $\calD_l$ with compact stabilizers isomorphic to $\SO(2)\times \O(19)$. This shows that  $\Gamma_l$ is a discrete subgroup of the automorphism group of $\calD_l$ (i.e. its stabilizers are finite subgroups). Thus the orbit space $\Gamma_l\backslash \calD_l$ has a structure of a complex variety (with quotient singularities). In fact, it has  a structure of a quasi-projective algebraic variety. Note that, although the period space $\calD_l$  is not connected, the orbit space is an irreducible algebraic variety. The reason is that the group $\Gamma_{l}$ contains an element which switches the two connected components. In fact, it is easy to see that there is an isomorphism of lattices
$$l^\perp \cong U^2\oplus E_8^2\oplus <-2d>,$$
where $<m>$ denote a lattice of rank 1 defined by the matrix $(m)$.  The isometry $\sigma$ of $l^\perp$ defined by being the minus-identity on one copy of $U$ and the identity on all other direct summands belongs to $\Gamma_l$. Since it switches the orientation of a positive definite 2-plane spanned by a vector $x\in U$ and a vector $y$ from another copy of $U$ such that $(x,x)> 0, (y,y) > 0$, the isometry $\sigma$ switches the two connected components of $\calD_l$.  

We see that the isomorphism class of $(X,h)$ defines a point 
$$p_{X,h}\in \Gamma_l\backslash \calD_l$$
and the Global Torelli Theorem from above asserts that $p(X,h) = p(X',h')$ if and only if $(X,h)$ and $(X',h')$ are isomorphic polarized surfaces. This is truly one of the deepest results in mathematics generalizing the Torelli theorem for  curves and polarized abelian varieties.

The Global Torelli Theorem can be restated in terms of the period map for a family $\pi:\calX \to T$ of polarized K3 surfaces, where we always assume that $T$ is connected. The cohomology groups of fibres $H^2(X_t,\bbZ)$ form a sheaf $\calH^2(\pi)$ of $\bbZ$-modules (the sheaf $R^2\pi_*(\bbZ_Y)$). A marking of the family is an isomorphism of sheaves $\phi:(L_{K3})_T\to \calH^2(\pi)$, where $(L_{K3})_T$ is the constant sheaf with fibre $L_{K3}$. For any $t\in T$ the map of fibres defines a marking $\phi_t$ of $X_t$. The  polarization classes $h_t$ of fibres define a section $h$ of $\calH^2(\pi)$ and we require that $\phi^{-1}(h) = l$. For any $t\in T$, the map of fibres $\phi_\bbC:
(L_{K3})_\bbC \to H^2(X_t,\bbC)$ defines the period point $p_{X_t,\phi_t}\in \calD_l$. This gives a map\begin{equation}\label{ppm}
p_{\pi,\phi}:T\to \calD_l.
\end{equation}
By Griffiths' theorem it is a holomorphic map of complex varieties. In general, a family does not admit a marking because the sheaf $\calH^2(\pi)$ need not be constant. However, it is isomorphic to a constant sheaf when the base $T$ is simply connected. So replacing $T$ by the universal covering $\tilde{T}$ we can define the period map \eqref{ppm}. Now for any two points $\tilde{t}$ and $\tilde{t'}$ in $\tilde{T}$ over the same point in $T$, the fibres of the family $\tilde{f}:Y\times_T\tilde{T}\to \tilde{T}$ over these points are isomorphic polarized surfaces. Therefore their images in $\calD_l$ lie in the same orbit of $\Gamma_l$. This defines a holomorphic map
$$p_\pi: T\to \Gamma_l\backslash \calD_l$$
called the {\it period map of a family}.  The Global Torelli Theorem asserts that the points in the same fibre of the period map are isomorphic polarized K3 surfaces. 

One can construct a coarse moduli space of polarized K3 surfaces. Let $X$ be an algebraic K3 surface which contains a Hodge class $h$ with $h^2 = 2d$. One can show that a line bundle $L$ with $c_1(L) = 3h$ defines an isomorphism from $X$ to a closed subvariety of $\bbP^{9d+1}$. Two different embeddings defined in this way  differ by a projective transformation of $\bbP^{9d+1}$. A construction from algebraic geometry based on the notion of a {\it Hilbert scheme} allows one to parametrize all closed subvarieties of 	$\bbP^{9d+1}$ isomorphic to an embedded K3 surface as above by some quasi-projective algebraic variety $H_{2d}$. The orbit space $H_{2d}/\text{PGL}(9d+2,\bbC)$ exists and represents a coarse moduli space of K3 surfaces (see \cite{SPS}). We denote it by $\calM_{2d}$. Then we have the period map 
\begin{equation}\label{ps}
\per:\calM_{2d} \to \Gamma_{l}\backslash \calD_l\end{equation}
The Global Torelli Theorem implies that this map is injective. A further result, the {\it Local Torelli Theorem} asserts that this map is a local isomorphism of complex varieties, and hence the period map is an isomorphism with its image. 

Before we state a result describing the image of the period map let us remind some terminology from algebraic geometry. Let $\Pic(X)$ be the Picard group of $X$. It is a subgroup of $H^2(X,\bbZ)$ spanned by the fundamental classes of irreducible 1-dimensional complex subvarieties (curves) on $X$. An integral  linear combination $\sum n_i [\Gamma_i]$ of such classes is cohomologous to zero if and only if the divisor $D = \sum n_i\Gamma_i$ is linearly equivalent to zero (i.e. the divisor of a meromorphic function). So $\Pic(X)$ can be identified with the group of divisor classes on $X$. Via the well-known relation between divisor classe and isomorphism classes of algebraic line bundles over $X$ (defined by the first Chern class map), we see another description of $\Pic(X)$ as the group of isomorphism classes of line bunldes, or invertible sheaves of $\calO_X$-modules. We identify $\Pic(X)$ with a sublattice $S_X$ of the unimodular lattice $H^2(X,\bbZ)$ with respect to the cup-product (the {\it Picard lattice}). 

Recall that a line bundle $L$ is called {\it ample} if sections of some positive tensor power $L^{\otimes m}$ define an embedding of $X$ in projective space. In this case the line bundle becomes isomorphic to the restriction of the line bundle corresponding to a hyperplane. A line bundle $L$ is ample if and only if its first Chern class $c_1(L)$ is the cohomology class of some Hodge form. A line bundle $L$ is ample if and only if 
$$c_1(L)^2 > 0, \ (c_1(L), [C]) > 0 ,$$ for any curve $C$.  We say that $L$ is {\it pseudo-ample} if the first condition holds but in the second one the equality is possible. It is known that the restriction of the cup-product $Q$ to $\Pic(X)$ has signature $(1,r-1)$, where $r = \rank~\Pic(X)$. This is because, by Lefschetz Theorem 
\begin{equation}\label{lef}
\Pic(X) = H^{11}(X)\cap H^2(X,\bbZ).\end{equation}
By the property of signature of the intersection form on $H^{11}$, any irreducible curve $C$ with $(c_1(L),[C]) = 0$ must satisfy $[C]^2 < 0$.  Since $K_X = 0$, the adjunction formula implies that $C\cong \bbP^1$ and $C^2 = -2$ (we will omit writing the brackets $[C]$ indicating the divisor class). Such a curve is called a {\it $(-2)$-curve}. Thus if $X$ has no $(-2)$-curves, any divisor class $D$ with $D\cdot D > 0$ is ample. It follows from the Riemann-Roch Theorem on algebraic surfaces that any divisor class $D$ with $D^2 = -2$ is either effective, i.e. the class of some curve, maybe reducible, or $-D$ is effective. 
It is known that a  line bundle $L$ is pseudo-ample if and only if sections of some positive tensor power of $L$ define a map from $X$ to a surface in $X'$ in some projective space such that all $(-2)$-curves $R$ with $(R,c_1(L)) = 0$ are blown down to singular points of $X'$ and ouside of the union of such curves the map is an isomorphism.

The next result describes the image of the period map \eqref{ps}. Let $\delta\in L_{K3}$ satisfy $(\delta,l) = 0, \delta^2 = -2$. Let $\calH_\delta$ denote the intersection of $\calD_l$ with the hyperplane 
$\{\bbC x\in \bbP((l^\perp)_\bbC) :(x,\delta) = 0\}.$ Let 
$$\Delta_l = \bigcup_{\delta^2 =-2,(\delta,l)= 0}\calH_\delta.$$
It is a union of countably many closed hypersurfaces in $\calD_l$. We call it the {\it discriminant locus} of $\calD_l$. Note that for any marked polarized K3 surface $(X,\phi)$ of degree $2d$ its period $p_{X,\phi}$ does not belong to any hyperplane $\calH_\delta$. In fact, otherwise 
$0 = (\delta,p_{X,\phi}) = (\phi(\delta),H^{20}(X))$ 
implies, by \eqref{lef} that $\phi(\delta)\in \Pic(X)$ and $\phi(\delta)$ or its negative is an effective divisor class. However 
$0 = (\delta,l) = (\phi(\delta),h)$ contradicts the assumption that $h = c_1(L)$ for some ample line bundle. So we see that the image of the period map is contained in the complement of the discriminant locus. A deep result, called the {\it Surjectivity Theorem} for the periods of polarized K3 surfaces asserts that the image is exactly the complement of $\Delta_l$. In fact, it gives more.

By weaking condition on the polarization one can consider pseudo-polarized K3 surfaces as  pairs $(X,h)$, where $h$ is a primitive divisor class corresponding to a pseudo-ample line bundle. One can extend the notion of a marked pseudo-polarized K3 surface $(X,\phi)$ and define its period point $p_{X,\phi}$  and extend the definition of the period map to families of {\it pseudo-polarized} K3 surfaces.  Then one proves the following.

\begin{thm} Any point in $\calD_l$ is realized as the period point $p_{X,\phi}$ of a marked pseudo-polarized K3 surface $(X,h)$. The image of the set $\{\delta:p_{X,\phi}\in \calH_\delta\}$ under the map $\phi$ is equal  to the set $\Delta(X,h)$ of divisor classes of curves $R$ such that $(R,h) = 0$. In particular, the periods of marked polarized K3 surfaces is equal to the set $\calD_l\setminus \Delta_l$.  The orbit space $\Gamma_l\backslash\calD_l$ is a coarse moduli space for pseudo-polarized K3 surfaces of degree $h^2=2d$. 
\end{thm}

Note that the set $\Delta(X,h)$ is finite because it is contained in the negative definite lattice equal to the orthogonal complement of $h$ in $S_X$. This shows that the set of hyperplanes $\calH_\delta$ is locally finite, i.e. each point lie in the intersection of only finite many hyperplanes. Let $\calR(X,h)$ be the sublattice of $S_X$ generated by the subset of  $\Delta(X,h)$. It is negaive definite (because it is a sublattice of the hyperbolic lattice $S_X$ orthogonal to an element $h$ with $h^2 > 0$) and is spanned by divisor classes $R$ with $R^2 = -2$. 
Let $[R]\in \Delta(X,h)$, replacing $R$  with $-R$ we may assume that $R$ is a curve on $X$. Writing it as a sum of its irreducible components $R_i$ we see that $(R_i,h) = 0$, so $R_i\in \Delta(X,h)$.  This shows that $\Delta(X,h)$ is spanned by the classes of $(-2)$-curves from $\Delta(X,h)$.  
All negative definite lattices generated by elements $x$ with $x^2 = 2$ lattices had been classified by E. Cartan.  They  coincide (after multiplying the value of the quadratic form by $-1$) with root lattices of semi-simple Lie algebras of types $A_n,D_n,E_n$ and their direct sums and described by the corresponding Cartan matrices and their direct sums. Recall that the pseudo-polarization $h$ defines a birational model of $X$ as a surface $X'$ with singular points whose minimal resolution is isomorphic to $X$. The singular points are rational double points of types corresponding to the irreducible components of the lattice $\calR_X$. The simplest case is when $\Delta(X,h)$ consists of one element. In this case, $X'$ has one ordinary singularity (of type $A_1$).  All possible lattices $\calR(X,h)$ which may occur have been classified only in the cases $h^2 = 2$ or $4$.

\section{Lattice polarized K3 surfaces}
Let $M$ be any abstract even lattice of signature $(1,r-1)$ (if $r > 1$ it is called a {\it hyperbolic lattice}). It is known  that the cone  
\begin{equation}
\label{cone}
 \calV_M = \{x\in M_\bbR: x^2 \ge 0\}
 \end{equation}
after deleting the zero vector consists of two connected components. In coordinates, the cone is linearly isomorphic to 
$$\{x = (x_1,\ldots,x_r)\in \bbR^r: x_1^2-x_2^2-\ldots-x_r^2 \ge 0\}.$$
The components are distinguished by the sign of $x_1$. 

Let  $\calV_M^0$ be one of the two ``halves'' of the cone $\calV_M$, and set 
$$M_{-2} = \{x\in M:x^2 = -2\}.$$ 
 Fix  a connected component $\calC(M)^+$ of 
 $$\calV_M^0\setminus \bigcup_{\delta\in M_{-2}}\delta^\perp.$$  
Let 
$$M_{-2}^+ = \{\delta\in M_{-2}: (m,\delta) > 0, \forall m\in \calC(M)^+\}.$$
This gives a decomposition
$$M_{-2} = M_{-2}^+\coprod M_{-2}^-,$$
where $M_{-2}^- = \{-\delta:\delta\in M_{-2}^+\}.$ We have
$$\calC(M)^+ = \{x\in \calV^0_M : (x,\delta) > 0, \forall \delta\in M_{-2}^+\}.$$

Take the Picard lattice $S_X$ as $M$.  Then we have a canonical decomposition of 
$$(S_X)_{-2} = (S_X)_{-2}^+ \coprod (S_X)_{-2}^-$$
where $(S_X)_{-2}^+$ consists of the classes of effective  divisors.  The corresponding connected component
is denoted by $\calK_X^0$ and is called {\it ample cone}.  We denote by $\calK_X$ the closure of $\calK^0_X$
in $\calV_{S_X}^0$.

 A {\it $M$-polarized  K3 surface} is a pair $(X,j)$ 
consisting of an algebraic K3 surface $X$ and a primitive
embedding of lattices 
$j: M\to \Pic(X)$ (where primitive means that the cokernel is a free abelian 
group) such that $j(\calC(M)^+)\cap \calK_X \ne \emptyset$. A $M$-marking is called 
{\it ample} if $j(\calC(M)^+)\cap \calK_X^0 \ne \emptyset$ (see \cite{Do}).

We say that two $M$-polarized K3 surfaces $(X,j)$ and $(X',j')$ are 
isomorphic if there exists an isomorphism of algebraic surfaces $f:X\to X'$  such 
that $j = f^*\circ j'$.

When $M$ is of rank 1 and is generated by an element $x$ with $x^2 = 2d$,  an $M$-polarized  (resp. ample $M$-polarized) K3 surface is a pseudo-ample (resp. ample) polarized K3 surface of degree $2d$.

 It turns out that the isomorphism classes of $M$-polarized K3 surfaces can be parametrized by a quasi-projective algebraic variety $\calM_M$. The construction is based on the period mapping.

First we fix a primitive embedding of $M$ in $L_{K_3}$. If it does 
not exist, $\calM_M = \emptyset$. We shall identify $M$ with its 
image. Then we consider a {\it marking} of $(X,j)$. It is an isomorphism of lattices
$\phi:L_{K3} \to H^2(X,\bbZ)$ such that $\phi|M = j$. 
An isomorphism of marked $M$-polarized surfaces is an isomorphism of surfaces $f:X\to X'$ such that 
$f^*\circ \phi' = \phi$.

Let $N = M^\perp$ be the orthogonal complement of $M$ in $L_{K3}$. Set
\begin{equation}
\calD_M = \{\bbC v\in \bbP(N_\bbC): v^2 = 0, (v,\bar{v}) > 0\}.
\end{equation}
In the case $M = \bbZ l$ is generated by one element, we see that $D_M\cong \calD_l$. Since $M$ is of signature $(1,r-1)$, the signature of the lattice $N$ is $(2,20-r)$. As before we see that 
$$\calD_M \cong G(2,N_\bbR)^+ \cong \O(2,20-r)/\SO(2)\times\O(20-r).$$
It is a disjoint union of two copies of a hermitian symmetric space of non-compact type IV of dimension $20-r$. 

Let $j(M) \subset S_X$, its orthogonal complement $j(M)^\perp$ in $H^2(X,\bbZ)$ is equal to $\phi(N)$ and does not depend on the choice of $\phi$. Set 
\begin{equation}\label{hsn}
V^{20} = H^{20}(X), \quad V^{11} = H^{11}(X)\cap \phi(N)_\bbC, \quad V^{02} = H^{02}(X).\end{equation}
This defines an AHS of weight 2 on the vector space $V= \phi(N)_\bbR$ with the period flag of type $\bff = (1,21-r,22-r)$. The polarization form is obtained from the  symmetric bilinear form defining the lattice structure on $N$. The space $\calD_M$ is isomorphic to the period space of such abstract Hodge structures. 

Since $\phi(M) \subset \Pic(X) \subset H^{11}(X)$, $\phi^{-1}(H^{20}(X))\subset N_\bbC$. This shows that
 $$p_{(X,\phi)}= \phi^{-1}(H^{20}(X))\in \calD_M.$$

Let 
$$N_{-2} = \{\delta\in N: \delta^2 = -2\}.$$
For any $\delta\in N_{-2}$  let 
$\delta^\perp$ denote   the 
hyperplane of vectors in $\bbP(N_\bbC)$ orthogonal to $\delta$. 
Let 
$\calH_\delta = \delta^\perp\cap \calD_M,$ and define 
$$\Delta_M =  \bigcup_{\delta\in N_{-2}}\calH_\delta.$$

\begin{thm}\label{surj} Let  $(X,j)$ be a  $M$-polarized K3 surface. Then 
\begin{itemize} 
\item[(i)] the $M$-polarization is ample if and only if for any choice of a marking 
$p_{X,\phi}\in \calD_{M}\setminus \Delta_M.$
\item[(ii)] If $(X,\phi)$ and $(X',\phi')$ are two  marked $M$-polarized K3 surfaces. Then $p_{X,\phi} = 
p_{X',\phi'}$ if and only if there exists an isomorphism of marked $M$-polarized surfaces $(X,\phi)\cong (X',\phi'\circ \alpha)$, where $\alpha$ is a product of reflections with respect to vectors $\delta\in N_{-2}$ such that $p_{X,\phi}\in  \calH_\delta$. 
\item[(iii)] Any point in $\calD_M$ is realized as the point $p_{X,\phi}$ for some $M$-polarized K3 surface.
\end{itemize}
\end{thm}

\begin{proof} Suppose $(X,j)$ is an ample $M$-polarized K3 surface and $\ell = p_{X,\phi}\in \calH_\delta$ for some $\delta\in N_{-2}$. This implies that $(\phi(\delta), H^{20}(X)) = 0$ and hence  $\phi(\delta)$ is the class of a divisor $R$ (maybe reducible) with $R^2 = -2$. Replacing $R$ with $-R$ if needed, we may assume that $R$ is a curve. Since $j(M)$ contains a polarization class and $[R]\in j(M)^\perp$ we get a contradiction. 

Conversely, suppose $p_{X,\phi}\not\in \Delta_M$ but $(X,j)$ is not ample. This means that $j(\calC(M)^+)\cap \calK_X\ne \emptyset$ but $j(\calC(M)^+)\cap \calK_X^0 =  \emptyset$. Equivalently, $j(\calC(M)^+)$ is contained in the boundary of $\calK_X$ defined by the equations $(x,[R]) = 0$,  where $R$ is a $(-2)$-curve. Since $\calC(M)^+$ spans $M$, this implies that there exists a $(-2)$-curve $R$ such that $\phi^{-1}([R])\in N_{-2}\cap \phi^{-1}(H^{20}(X)^\perp)$. Therefore $p_{X,\phi} = \phi^{-1}(H^{20}(X))\in \calH_\delta, $where $\delta = \phi^{-1}([R])$, contradicting the assumption. 

The last two assertions follow from the Global Torelli Theorem of Burns-Rappoport for non-polarized K\"ahler surfaces (for a good exposition of their work see \cite{Ge}).  We refer for the proof of these assertions to \cite{Do}. Note that the fixing $\calC(M)^+$ in the definition of a $M$-polarized surface is used for the proof of assertion (ii).  
\end{proof}

Now let us see how the period point $p_{X,\phi}$ depends on the 
marking $\phi$. For any isometry $\sigma$ of the lattice $L_{K3}$ which
is identical on $M$ we obtain another marking $\sigma\circ \phi$ of
the $M$-polarized surface $(X,j)$. Let $\Gamma_M$ denote the group of
such isometries. We have a natural injective homomorphism $\Gamma_M\to
\O(N)$ and one can easily determine the image. Let $D(N) = N^*/N$ be
the discriminant group of $N$. The group $\O(N)$ acts naturally on
$N^*$ leaving $N$ invariant. This defines a homomorphism $\O(N)\to
\Aut(D(N))$ and we have
$$\Gamma_M \cong \O(N)^* := \Ker \bigl(\O(N)\to \Aut(D(N))\bigr).$$
This is because every isometry from the right-hand side can be extended to 
an isometry of $L$ acting identically on $M = N^\perp.$ The group
$\Gamma_M$ is a discrete subgroup of the group of automorphims of 
$\calD_M$ and the orbit space $\Gamma_M\backslash \calD_M$ is a quasi-projective algebraic variety.  The image of $p_{X,\phi}$ in the orbit space $\Gamma_M\backslash\calD_M$ depends only
on $(X,j)$. It is denoted by $p_{X,j}$ and is called the {\it period of 
$M$-polarized surface $X$}.

Let $\pi:Y\to T$ be a family of K3 surfaces parametrized by some algebraic variety 
$T$. The cohomololy groups $H^2(X_t,\bbZ), t\in T,$ form a locall coefficient
system $\calH^2$ on $T$. Assume that there exists an embedding $j:M_T\hookrightarrow 
\calH^2$ of the constant coefficient system $M_T$ such that for each $t\in T$ the
embedding
$j_t:M\to H^2(X_t,\bbZ)$ defines an M-polarization of $X_t$. We say then that
 $\pi$ 
is a {\it family of $M$-polarized $K3$ surfaces}. We leave to the reader to extend the notion of the period map of polarized K3 surfaces to the case of $M$-polarized K3 surfaces
$$p_\pi:T\to \Gamma_M\backslash\calD_M.$$

This is called the {\it period map} of the family $\pi:\calX\to T$.

Note that the group $\Gamma_M$ contains reflections $s_\delta:x\mapsto (x,\delta)\delta$ for any $\delta\in N_{-2}$. Applying Theorem \ref{surj}, we obtain that for any  $M$-polarized K3 surfaces $(X,j)$ and $(X',j')$
\begin{equation}\label{surj2}
p(X,j) = p(X',j') \Longleftrightarrow (X,j) \cong (X',j').
\end{equation} 

This implies that $\Gamma_M\backslash \calD_{M}$ (resp. $\Gamma_M\backslash \calD_{M}^0$) can be taken as  a coarse moduli space $\calM_{M,K3}$ (resp. $\calM_{M,K3}^a$) of   $M$-polarized (resp. ample $M$-polarized) K3 surfaces.

Notice that there is only finitely many orbits of $\Gamma_M$ on the set 
$N_{-2}.$
So, 
$$\Gamma_M\backslash \Delta_M
$$ is a closed subset in Zariski topology of 
 $\Gamma_M\backslash\calD_M$.
We know that in the case when $M$ is of rank 1, the orbit space is irreducible because $\Gamma_M$ contains an element swithching two connected components of $\calD_M$. The same is true (and the same proof works)  if the lattice $N$ contains a direct summand isomorphic to $U$ \cite{Do}.

Let $(X,j)$ be a  $M$-polarized K3 surface. Note that in general $j(M) \ne S_X$. As we explained, a surface $(X,j)$ is ample $M$-polarized if and only if $j(M)^\perp \cap S_X$ contains no curves  $R$ with $R^2 = -2$.   As in the case of a pseudo-polarized K3 surface we can define the {\it degeneracy lattice} $\calR(X,j)$  as the span of $(-2)$-curves  in $j(M)^\perp$.

\section{Eigenperiods of algebraic K3 surfaces}\label{ab}

First we racall the fundamental result due to Nikulin \cite{N2}.

\begin{thm}
Let $X$ be an algebraic K3 surface and let $\omega_X$ be a nowhere vanishing holomorphic 2-form on $X$.  Let
$$\alpha : {\rm Aut}(X) \to \bbC^*$$ 
be a homomorphism defined by $g^*(\omega_X) = \alpha(g)\cdot \omega_X$ for $g \in {\rm Aut}(X)$.
Then

$(1)$  the image of $\alpha$ is a cyclic group of a finite order $m$.

$(2)$  Assume that $\alpha(g)$ is a primitive $m$-th root of 1.  Then $g^*$ has no non-zero fixed vectors in $T_X\otimes \bbQ$. 
\end{thm}

Let $C_m$ be a cyclic group of order $m$ of automorphims of  a $M$-polarized K3 surface $(X,j)$. 
We assume that the restriction of $g^*$ to $j(M)$ is the identity map  for any $g\in C_m$.
Fixing a marking $\phi:L_{K3}\to H^2(X,\bbZ)$ such that $\phi|M = j$ we obtain a  homomorphism
$$\rho_\phi:C_m\to \O(L_{K3}), \quad g\mapsto  \phi^{-1}\circ g^{-1}{}^*\circ \phi.$$ 

Fix a homomorphism  $\rho:C_m\to \O(L_{K3})$ such that 
\begin{equation}\label{fix}
M = L_{K3}^\rho:= \{x\in L_{K3}: \rho(a)(x) = x, \forall a\in C_m\}.
\end{equation}
Define a $(\rho,M)$-marking of $(X,j)$ as a marking $\phi:L_{K3}\to H^2(X,\bbZ)$ with $\phi|M = j$  such that $\rho = \rho_\phi$.

Let $\chi\in \hat{C}_m$ be the unique character of $C_m$ such that 
$$H^{20}(X) \subset H^2(X,\bbC)(\chi).$$
Then 
$$p_{X,\phi} = \phi^{-1}(H^{20}(X))\subset  N_\bbC(\chi)$$
where $N = M^{\perp}$.

Assume that $\chi$ is not a real character. Then $N_\bbC(\chi)$ is an isotropic subspace of $N_\bbC$ with respect to the quadratic form $Q$ defined by the lattice $N$. The restriction of the hermitian form $Q(x,\bar{y})$ to $p_{X,\phi}$ is positive. Let 
$$\calD_{M}^{\rho,\chi} = \{\bbC x\in N_\bbC(\chi):Q(x,\bar{x}) > 0.\}.$$
We know that 
$$\calD_{M}^{\rho,\chi} \cong \bbB_{d(\chi)},$$
where $d(\chi) = \dim N_\bbC(\chi)-1.$

If $\chi$ is a real character, then $N_\bbC(\chi) = N_\bbR(\chi)\otimes\bbC$ and we set
$$\calD_{M}^{\rho,\chi} = \{\bbC x\in N_\bbC(\chi):Q(x,x) = 0, Q(x,\bar{x}) > 0\}.$$
This is nothing but a type IV Hermitian symmetric space of dimension $d(\chi)-1$.

There is a natural inclusion $\calD_{M}^{\rho,\chi}\subset \calD_M$ corresponding to the inclusion of $N_\bbC(\chi)\subset N_C$. 

The group $C_m$ acts on $\calD_M$ via the restriction of the  representation $\rho$ to the homomorphism 
$\rho':C_m\to \O(N)$. We denote by 
$\calD_M^\rho$ the set of points which are fixed under all $\rho(a), a\in C_m$. It is clear that 
$$\calD_M^\rho = \coprod_{\chi\in \hat{C}_m}\calD_{M}^{\rho,\chi}.$$

\begin{thm}\label{surj1} Let $(X,\phi)$ be a $(\rho,M)$-marked ample $M$-polarized K3 surface such that $H^{20}(X)\subset H^2(X,\bbC)(\chi).$ Then
 $p_{X,\phi}\in \calD_M^{\rho,\chi}\setminus \Delta_M$. Conversely, any point in  $\calD_M^{\rho,\chi}\setminus \Delta_M$ is equal to $p_{X,\phi}$ for some $(\rho,M)$-marked ample $M$-polarized K3 surface $(X,\phi).$ 
 \end{thm}
 
 \begin{proof} The first part  was already explained. Let $\ell$ be a point in $\calD_M^{\rho,\chi}\setminus \Delta_M$. It represents a line in $N_\bbC$ which is fixed under $\rho(a), a\in C_m$. Let $\phi:L_{K3} \to H^2(X,\bbZ)$ be an ample $M$-marking of $X$ such that  $\phi(\ell) = H^{20}(X)$.  Then $\phi^{-1}\circ \rho(a)\circ \phi$ acts identically on $j(M)$ and preserves $H^{20}(X)$.  Since $(X,j)$ is ample $M$-polarized, $j(M)$ contains an ample class which is preserved under $\rho(a)$.  By Global Torelli Theorem $\rho(a) = g^*$ for a unique $g\in \Aut(X)$.  This defines an action of $C_m$ on $X$ and a  $(\rho,M)$-marking of $X$. Obviously, $H^{20}(X)$ and $p_{X,\phi}$ are eigenspaces corresponding to the same character of $C_m$.
\end{proof}

Let 
\begin{equation}\label{gam}
\Gamma_{M,\rho} = \{\sigma\in \Gamma_M: \sigma\circ \rho(a) = \rho(a)\circ \sigma, \forall a\in C_m\}.
\end{equation}

Obviously, the group $\Gamma_{M,\rho}$ leaves the connected components $\calD_{M}^{\rho,\chi}$ of $\calD_M^\rho$ invariant. 
Applying Theorem \ref{surj1}, we obtain

\begin{thm}  The orbit space $\Gamma_{M,\rho}\backslash( \calD_{M}^{\rho,\chi}\setminus \Delta_M) $  parametrizes isomorphism classes of ample $(\rho,M)$-polarized  K3 surfaces such such that $H^{20}(X)\subset H^2(X,\bbC)(\chi)$. 
\end{thm}

\begin{rem}  
Let 
\begin{equation}\label{gam}
\tilde{\Gamma}_{M,\rho} = \{\sigma\in \O(N) : \sigma\circ \rho(a) = \rho(a)\circ \sigma, \forall a\in C_m\}.
\end{equation}
If the canonical map $\tilde{\Gamma}_{M,\rho} \to \O(D(N))$ is surjective, the quotient $\tilde{\Gamma}_{M,\rho}/\Gamma_{M,\rho}$ is
isomorphic to $\O(D(N)) \cong \O(D(M))$.
The orbit space $\tilde{\Gamma}_{M,\rho}\backslash( \calD_{M}^{\rho,\chi}\setminus \Delta_M) $  
(resp. $\Gamma_{M,\rho}\backslash( \calD_{M}^{\rho,\chi}\setminus \Delta_M) $)
sometimes parametrizes the moduli space of some varieties (resp. some varieties with a marking).
See  section \ref{example}.
\end{rem}

\begin{lem}
Let $p_{X,\phi}\in \calD_M^\rho\cap \Delta_M$.  Then $\rho = \rho_\phi$ is not represented by any
automorphism of $X$.
\end{lem}

\begin{proof}
Assume that $p_{X,\phi}\in \calD_M^\rho\cap \Delta_M$ and $\rho = \rho_\phi$ is a
representation of $C_m$ acting on $X$ as automorphisms.  
Since $C_m$ is finite, by averaging, we can find a polarization class $h$ such that $g^*(h) = h$ for all $g\in C_m$. By condition \eqref{fix}, $h\in j(M)$ which contradicts the assumption that $(X,j)$ is not ample. 
\end{proof}

In the following we study the locus $\calD_M^\rho\cap \Delta_M$
under the assumption that $m$ is an odd prime number $p$.  Recall that $C_p$ acts trivially on $M$ and acts effectively on 
$N = M^{\perp}$.  This implies that $D(M) \cong D(N)$ are $p$-elementary abelian groups.
Consider the  degeneracy lattice $\calR(X,j)$ of $(X,j)$. 
If $r \in \calR(X,j)$ is a $(-2)$-vector, then $\la r, \omega_X\ra = 0$ implies that $\la \rho^i(r), \omega_X\ra = 0$.
Hence for any marking $\phi$ of $(X,j)$ its pre-image 
$R = \phi^{-1}(\calR(X,j))$ is an invariant sublattice of $L_{K3}$ under $\rho$ orthogonal to $M$. 
We consider the generic case, that is, $R$ is generated by $r, \rho(r),..., \rho^{p-2}(r)$.

\begin{lem}\label{p-elem}
Let $r \in \calR(X,j)$ with $r^2=-2$.  Let $R$ be the lattice generated by $r, \rho(r),..., \rho^{p-2}(r)$.  Then
$R$ is $\rho$-invariant and isomorphic to the root lattice $A_{p-1}$.
Moreover $\rho$ acts trivially on $D(R)$.
\end{lem}

\begin{proof}  Since $R$ is generated by $(-2)$-vectors and negative definite, $R$ is a root lattice of rank $p-1$.
The $\rho$-invariance is obvious.
Let $m_i = \langle r, \rho^i(r) \rangle$.  
Then the definiteness implies that $|m_i| \leq 1$.
Obviously $m_1= m_{p-1}, m_2 = m_{p-2},..., m_{(p-1)/2} = m_{(p+1)/2}$. 
Since $-2 = r^2 = \langle r, -\sum_{i=1}^{p-1} \rho^i (r) \rangle = -2m_1- \cdot  \cdot  \cdot - 2m_{(p-1)/2}$, 
$m_i = 1$ for some $1\leq i\leq (p-1)/2$.  If $m_i = m_j= 1$ for $1\leq i\not=j\leq (p-1)/2$, 
then $r, \rho^i(r), \rho^j(r), \rho^{p-i}(r), \rho^{p-j}(r)$ generate
a degenerate lattice. This contradicts the definiteness of $R$.  Thus there exists a unique $i$ such that
$m_i=1$ and $m_j =0$ for $1\leq i\leq (p-1)/2$.  By changing $\rho$ by $\rho^i$, we may assume that $m_1=1$.
Then $\langle \rho^i(r), \rho^{i+1}(r) \rangle = 1$ and hence 
$r,\rho(r),..., \rho^{p-2}(r)$ generate a root lattice $A_{p-1}$.  The last assertion follows from the fact that
$\rho$ fixes a generator $(r + 2\rho(r) + \cdot \cdot \cdot + (p-1)\rho^{p-2}(r))/p$ of $D(A_{p-1})$.
\end{proof}

Let $R$ be as in Lemma \ref{p-elem}.
Let $M'$ be the smallest primitive sublattice of $L_{K3}$ containing $M\oplus R$. 
We have a chain of lattices
$$0 \subset M\oplus R \subset M'\subset M'{}^* \subset (M\oplus R)^* .$$
Since any $\rho(a)$ acts identically on $D(M\oplus R) = D(M)\oplus D(R)$ (Lemma \ref{p-elem}), 
it acts identically on $D(M')$. Let $N' = M'{}^\perp$. Since $N'$ is the orthogonal complement of $M'$ in a unimodular lattice, the discriminant groups $D(N')$ and $D(M')$ are isomorphic and the action of $C_m$ on $D(N')$ is also identical.   This implies that $(1_{M'}\oplus \rho(a)|N')$ can be extended to an isometry 
$\rho'(a)$ of $L_{K3}$. This defines a homomorphism $\rho':C_m\to \O(L_{K3})$ with the sublattice of fixed vectors $L_{K3}^{\rho'}$ equal to $M'$. Let $j' = \phi|M':M'\to \Pic(X)$ and let $M'_{-2} = M_{-2}'{}^+\coprod M_{-2}'{}^-$ be the  partition of $M_{-2}'$ which extends the partition of $M_{-2}$ and is defined by the condition that 
$j'(M_{-2}'{}^+)$ is contained in the set of effective divisor classes. Then $(X,j')$ is a $M'$-polarized K3 surface. The polarization is ample if $p(X,j')\not\in \calH_{\delta}$ for any $\delta\in N'_{-2}$.  If this condition holds, the pair $(X,j')$ is a $(\rho',M')$-polarized K3 surface.  By induction,
we can prove the following result.

\begin{thm}\label{deg}
Assume that $m$ is a prime number $p$.
For any point $x \in \calD_M^{\rho,\chi}$, let $R(x)$ denote the sublattice of $N$ generated by the vectors $\delta\in N_{-2}$ such that $x\in \calH_\delta$.  Then $R(x)$ is a $p$-elementary root lattice, that is, $R(x)^*/R(x)$ is a $p$-elementary
abelian group.
Let $M'$ be the smallest primitive sublattice of $L_{K3}$ generated by $M$ and $R(x)$. 
Then $x = p_{X,\phi}$ for some 
marked ample $(\rho',M')$-polarized $K3$ surface $X$, where $\rho'$ is a representation of $C_p$ on $L_{K3}$ with $L_{K3}^{\rho'} = M'$. The lattice $R(x)$ is isomorphic to the degeneracy lattice of the 
$M$-polarized $K3$ surface $(X,j)$, where $j = \phi|M$.
\end{thm}

\section{Examples}\label{example}

In this section, we shall give examples of eigenperiods of K3 surfaces.  For a lattice $L$ and an integer $m$, we
denote by $L(m)$ the lattice over the same $\bbZ$-module with the symmetric bilinear form
multiplied by $m$.

\begin{ex} Let $C_2 = (g)$ be a cyclic group of order 2. Write $L_{K3}$ as a direct sum of two lattices 
$M = E_8\oplus U$ and $N =  E_8\oplus U\oplus U$ and define $\rho:C_2\to \O(L_{K3})$ by 
$\rho(g)$ to be the identity on the first summand and the minus identity on the second one. Let $f:X\to \bbP^1$ be an elliptic K3 surface with a section $s$ (taking as the zero in the Mordell-Weil group of sections) and assume that one of the fibres is a reducible fibre of type II* in Kodaira's notation.  Choose an $M$-marking $\phi$ of $X$ by fixing an isomorphism from $E_8$ to the subgroup of $S_X$ generated by components of the reducible fibre not intersecting $s$ and an isomorphism from $U$ to the subgroup of $S_X$ to be the pre-image under $\phi$ generated by any fibre and the section $s$.   The partition of $M_{-2}$ is defined by taking $M_{-2}^+$ to be the pre-image under $\phi$ of the set of divisor classes of irreducible components of the reducible fibre and the section. 
Now the coarse moduli space  of isomorphism classes of $M$-polarized K3 surfaces consists of elliptic K3-surfaces with a  section and a fibre of type $II^*$. Each such surface admits an automorphism of order 2 defined by the inversion automorphism $x\mapsto -x$ on a general fibre with respect to the group law defined by the choice of a section.   The  $M$-polarization is ample if and only if  the elliptic fibration has only one reducible fibre.  Let $\chi\in \hat{C_2}$ be defined by $\chi(g) = -1$.  Then $\calD_M^{\rho,\chi} \cong \calD_M$ and $\Gamma_{M,\rho} = \Gamma_M$. An ample $M$-polarized K3 surface  is automatically $(\rho,M)$-polarized K3 surface. The degeneracy lattice $\calR(X,j)$ is the sublattice of $\Pic(X)$ generated by components of new reducible fibres not intersecting the section.  Using Theorem \ref{deg} we see that any $M$-polarized K3 surface admits a structure of a $(\rho',M')$-polarized K3 surface, where $M'$ is isomorphic to the sublattice of $S_X$ generated by components of fibres and a section and $\rho'$ acts identically on $M'$ and the minus-identity on its orthogonal complement.
\end{ex}

\begin{ex}
(6 points on $\bbP^1$) 
 Let $X'$ be a surface in weighted projective space $\bbP(1,1,2,2)$ given by an equation $f_6(x_0,x_1)+x_2^3+x_3^3 = 0$, where $f_6$ is a homogenous form of degree 6 without multiple roots.  The surface $X'$ has 3 ordinary nodes $(0,0,1,a)$, where $a^3 = -1$. A minimal resolution of these nodes is a K3 surface $X$ on which the cyclic group $C_3 = \mu_3$ of  3th roots of unity acts by 
$(x_0,x_1,x_2,x_3)\mapsto (x_0,x_1,\alpha x_2,\alpha x_3)$. Another way to define this surface is to consider a quadratic cone in $\bbP^3$, take its transversal intersection with a cubic surface, and then to define $X'$ to be the triple cover of the cone  branched along the intersection curve $C$. The pre-image of the ruling of the cone defines an elliptic fibration $f:X\to \bbP^1$. Since $C$ has 6 ramification points of order 3 over the roots of $f_6$, the fibration has 6 reducible fibres of Kodaira type $III$. The three exceptional curves of the resolution $X\to X'$ form a group of sections isomorphic to $\bbZ/3\bbZ$. An explicit computation of the Hodge structure on hypersurfaces in a weighted projective space (see \cite{Do2}) shows that 
$H^{20}(X) \subset H^2(X,\bbC)(\chi)$, where $\chi(\alpha) = \alpha^2$ and 
$$\dim H^{11}(X)(\chi) = 3, \ \dim H^{11}(X)^{C_3} = 12.$$
  Let $M$ be the primitive sublattice of $S_X$ generated by the classes of irreducible components of fibres, and  sections. Standard arguments from the theory of elliptic surfaces yield that 
$M \cong U\oplus E_6\oplus A_2^3$ and its orthogonal complement $N \cong A_2(-1)\oplus A_2^3.$
We see that 
$$\calD_M^{\rho,\chi} \cong \bbB_3.$$
The lattice $N$ has a natural structure of a module over the ring of Eisenstein numbers $\bbZ[e^{2\pi i/3}]$ equipped with a hermitian form over this ring. One shows that
$$\Gamma_{M,\rho} \cong \SU(4,\bbZ[e^{2\pi i/3}]).$$
Fix a primitive embedding of $M$ in $L_{K3}$ and choose $\calC(M)^+$ to be the nef cone in $\Pic(X_0)$ for a fixed $X$ as above. Let $\rho$ be the representation of $C_3$ on $L_{K3}$ which acts as the identity on $M = U\oplus E_6\oplus A_2^3$ and acts on $N$ as follows. Choose  a standard basis $(r_1,r_2)$ of $A_2$ corresponding to the vertices of the Dynkin diagram. We want to define an action $\rho$  of the group $C_3$  on $N$ such that $N^{C_3} = \{0\}$ and $C_3$ acts trivially on the discriminant group of $N$ (because it acts trivially on $M$). This easily implies that each direct summand must be invariant, and the action on it has no fixed vectors and the generator of $D(A_2)$ equal to $(r_1+2r_2)/3$ is invariant. It is easy to see that the representation on $A_2$ with this property must be isomorphic to the one given by
$\rho(e^{2\pi i/3})(r_1,r_2) = (r_2,-r_1-r_2)$. Now let $\rho(e^{2\pi i/3})$ act on $N$ as the direct sum of the previous action on each direct summand.  In this way we obtain that the coarse moduli space of $(M,\rho)$-polarized K3 surfaces is isomorphic to the ball quotient
$$\Gamma_{M,\rho}\backslash\calD_M^{\rho,\chi} \cong \SU(4,\bbZ[e^{2\pi i/3}])\backslash\bbB_3.$$
Note that, using the hypergeometric functions  with $\bldmu = (1/3,1/3,1/3,1/3,1/3,1/3)$  the same quotient is isomorphic to the space $\calP_{1,\bfk}$ of ordered 6 points on $\bbP^1$, where 
$\bfk = (1,1,1,1,1,1)$. It is also isomorphic to the moduli space of principally polarazied abelian varieties of dimension 4 with action of a cyclic group of order 3 of type $(1,3)$. 
\end{ex}

\begin{ex}
(del Pezzo surfaces of degree 2)
Let $R$ be a smooth del Pezzo surface of degree 2.
Recall that the anti-canonical model of a smooth del Pezzo surface of degree 2 is a double cover of
$\bbP^2$ branched along a smooth plane quartic curve $C$.  Let $f_4(x,y,z)$ be the homogeneous polynomial
of degree 4 defining $C$.  Let $X$ be a quartic surface in $\bbP^3$ defined by
$$t^4 = f_4(x,y,z).$$  Then $X$ is a K3 surface which is the 4-cyclic cover of $\bbP^2$ branched along $C$.
Obviously $X$ admits an automorphism $g$ of order 4.  If $C$ is generic, then the Picard lattice $S_X$ 
(resp. the transcendental lattice $T_X$) is isomorphic to $M = U(2)\oplus A_1^{\oplus 6}$ (resp.
$N=U(2)\oplus U(2) \oplus D_8\oplus A_1^{\oplus 2}$).  In this case, $(g^*)^2$ acts trivially on $M$, but $g^*$ 
does not.  
The above correspondence gives an isomorphism between the moduli space of smooth del Pezzo surfaces of
degree 2 and $(\bbB_6\setminus \Delta_M)/\tilde{\Gamma}_{M,\rho}$.  
The quotient $(\bbB_6\setminus \Delta_M)/\Gamma_{M,\rho}$ is isomorphic to the moduli space of
{\it marked} smooth del Pezzo surfaces of degree 2.  The natural map 
$$\Gamma_{M,\rho}\backslash \bbB_6 \to \tilde{\Gamma}_{M,\rho}\backslash\bbB_6$$
is a Galois covering with $\O(D(M)) \cong \Sp(6,\bbF_2)$ as its Galois group.
We remark that this correspondence gives a uniformization of the moduli space of
non-hyperelliptic curves of genus 3.  For more details we refer the reader to the paper \cite{K1}.
\end{ex}

\begin{ex}
(8 points on $\bbP^1$)
This case is a degenerate case of the previous example since 8 points on $\bbP^1$ corresponds to a hyperelliptic
curve of genus 3.
Let $\{ (\lambda_{i} : 1)\}$ be a set of distinct 8 points on the projective line.
Let $(x_{0}: x_{1}, y_{0}: y_{1})$ be the  bi-homogenious coordinates on $\bbP^1 \times \bbP^1$.  
Consider a smooth divisor $C$ in $\bbP^1 \times \bbP^1$ of bidegree $(4,2)$ given by
\begin{equation}\label{embed}
y_{0}^{2} \cdot \prod^{4}_{i=1} (x_{0} - 
\lambda_{i} x_{1}) + y_{1}^{2} \cdot 
\prod^{8}_{i=5} (x_{0} - \lambda_{i} x_{1}) = 0.
\end{equation}
Let $L_{0}$ (resp. $L_{1}$) be the divisor defined
by $y_{0} = 0$ (resp. $y_{1} = 0$).
Let $\iota$ be an involution of 
$\bbP^1 \times \bbP^1$ given by
\begin{equation}\label{invol}
(x_{0}: x_{1}, y_{0}: y_{1})  \longrightarrow
(x_{0}: x_{1}, y_{0}: -y_{1})
\end{equation}
which preserves $C$ and $L_{0}, L_{1}$.  
Note that the double cover of $\bbP^1 \times \bbP^1$ 
branched along $C + L_{0} +L_{1}$ has 8 rational double points of type
$A_{1}$ and its minimal resolution $X$ is a $K3$ surface.   
The involution $\iota$ lifts to an automorphism $\sigma$ of order 4.
The projection
$$(x_{0}: x_{1}, y_{0}: y_{1})  \longrightarrow (x_{0}: x_{1})$$
from $\bbP^1 \times \bbP^1$ to $\bbP^1$ induces an
elliptic fibration 
$$\pi : X \longrightarrow \bbP^1$$
which has 8 singular fibers of type $III$ and two sections.
In this case, $M \simeq U(2) \oplus D_{4} \oplus D_{4}$ and
$N \simeq U \oplus U(2) \oplus D_{4} \oplus D_{4}.$
Thus we have an isomorphism between the moduli space of {\it ordered} (distinct) 8 points on $\bbP^1$ and
$\Gamma_{M,\rho}\backslash (\bbB_5\setminus \Delta)$.  
This ball quotient appeared in Deligne-Mostow's list \cite{DM}.
The group $\O(D(M)) \cong S_8$ naturally acts on the both spaces.  Taking the quotient by $S_8$, we have an
isomorphism between the moduli space of (distinct) 8 points on $\bbP^1$ and
$\tilde{\Gamma}_{M,\rho}\backslash (\bbB_5\setminus \Delta)$.
For more details we refer the reader to the paper \cite{K3}.
\end{ex}

\begin{ex}
(del Pezzo surfaces of degree 3)  Let $S$ be a smooth del Pezzo surface.  The anti-canonical model
of $S$ is a smooth cubic surface.  Let $l$ be a line on $S$.  Consider the conic bundle on $S$ defined by
hyperplanes through $l$.  The line $l$ is a 2-section of this pencil.  For generic $S$, this pencil has 5
degenerate members.  Moreover there are two fibers tangent to $l$.  These give homogeneous polynomials 
$f_5(x,y)$, $f_2(x,y)$ of degree 5 and 2 in two variables.
Now consider the plane sextic curve $C$ defined by
$$z(f_5(x,y) + z^3f_2(x,y)) = 0.$$
Take the double cover of $\bbP^2$ branched along $C$ whose minimal resolution $X$ is a K3 surface.
The multiplication of $z$ by $e^{2\pi \sqrt{-1}/3}$ induces an automorphism $g$ of $X$ of order 3.
In this case, $M \simeq U \oplus A_2^{\oplus 5}$ and
$N \simeq A_2(-1) \oplus A_2^{\oplus 4}.$
The non trivial fact is that the K3 surface $X$ is independent of the choice of a line $l$.
Thus we have an isomorphism between the moduli space of smooth cubic surfaces and
$\tilde{\Gamma}_{M,\rho}\backslash (\bbB_4\setminus \Delta)$.  
The group $\O(D(M)) \cong W(E_6)$ appears as the Galois group of the covering 
$$\Gamma_{M,\rho}\backslash \bbB_4 \to \tilde{\Gamma}_{M,\rho}\backslash\bbB_4.$$
We remark that the pencil of conics on $S$ induces an elliptic fibration on $X$ with five singular fibers of type IV and
two singular fibers of type II.  This gives a set of 7 points on $\bbP^1$.  The ball quotient 
$\Gamma_{M,\rho}\backslash \bbB_4$ appeared in Deligne-Mostow's list \cite{DM}.  
For more details we refer the reader to the paper \cite{DGK}. Note that another approach (see \cite{ACT}) consists of attaching to a cubic surface with equation $f(x_0,x_1,x_2,x_3) = 0$ the cubic hypersurface with equation $f(x_0,x_1,x_2,x_3)+x_4^3 = 0$ and considering its intermediate Jacobian variety. It is a principally polarized abelian variety of dimension 5 with cyclic group of order 3 acting on it with type $(1,4)$. As we know from section \ref{sec6}, such varieties are parametrized by a 4-dimensional ball.

\end{ex}

\begin{ex}
(del Pezzo surfaces of degree 4)
Let $S$ be a smooth del Pezzo surface of degree 4.  Its anti-canonical model is the complete intersection of
two quadrics $Q_1, Q_2$ in $\bbP^4$.  It is known that $Q_1$ and $Q_2$ can be diagonalized simultaneously,
that is,  we may assume 
$$Q_1= \{\sum z_i^2 = 0 \}, \\ Q_2 = \{\sum \lambda_i z_i^2 = 0\}.$$
The discriminant of the pencil of quadrics $\{ t_1Q_1 + t_2Q_2 \}_{(t_1:t_2)}$ is distinct 5
points $\{(\lambda_i : 1)\}$ on $\bbP^1$.  Conversely distinct 5 points on  $\bbP^1$ gives the intersection of two
quadrics.  Thus the moduli space of smooth del Pezzo surfaces of degree 4 is isomorphic to the moduli of
distinct 5 points on $\bbP^1$.  Next we construct a K3 surface from distinct 5 points $\{(\lambda_i : 1)\}$ on $\bbP^1$. 
Let $C$ be a plane sextic curve defined by
$$x_{0}^{6} = x_0 \prod_{i=1}^5 (x_1 - \lambda_i x_2).$$
Obviously $C$ is invariant under a projective transformation of order 5 by the multiplication of $x_0$ by a primitive
5-th root of 1.
Take the double cover of $\bbP^2$ branched along $C$ whose minimal resolution $X$ is a K3 surface.
The $X$ has an automorphism of order 5 induced from the above projective transformation.
In this case $M$ is of rank 10 and $D(M) = (\bbZ/5\bbZ)^3$, and
$$N \cong 
\begin{pmatrix}0&1
\\1&0
\end{pmatrix}
\oplus
\begin{pmatrix}2&1
\\1&-2
\end{pmatrix}
\oplus A_4 \oplus A_4.
$$
Thus we have an isomorphism between the moduli space of {\it ordered} (distinct) 5 points on $\bbP^1$ and 
$\Gamma_{M,\rho}\backslash (\bbB_2\setminus \Delta)$.  The ball quotient 
$\Gamma_{M,\rho}\backslash (\bbB_2\setminus \Delta)$ appeared in Deligne-Mostow \cite{DM}.
The group $\O(D(M)) \cong S_5$ naturally acts on the both spaces and the quotients give an isomorphism
between the moduli space of smooth del Pezzo surfaces of degree 4 and
$\tilde{\Gamma}_{M,\rho}\backslash (\bbB_2\setminus \Delta)$.
For more details we refer the reader to the paper \cite{K4}.
\end{ex}

\begin{ex}
(Curves of genus 4)
Let $C$ be a non-hyperelliptic curve of genus 4.  The canonical model of $C$ is the intersection of
a smooth quadric $Q$ and a cubic surface $S$ in $\bbP^3$.  By taking the triple cover of $Q$ branched along
$C$, we have a smooth K3 surface $X$ with an automorphism of order 3.  In this case $M = U(3)$ and
$N = U\oplus U(3) \oplus E_8 \oplus E_8$.  
Thus we have an isomorphism between the moduli space of non-hyperelliptic curves of genus 4 and 
$\tilde{\Gamma}_{M,\rho}\backslash (\bbB_9\setminus \Delta)$.  We remark that a ruling of $Q$ induces an elliptic fibration on $X$
with 12 singular fibers of type II.  This gives a set of 12 points on $\bbP^1$.
The arithmetic subgroup  
$\Gamma_{M,\rho}$ is isogenous to the complex reflection group associated to 12 points on $\bbP^1$ 
appeared in Mostow \cite{Mo}.
For more details we refer the reader to the paper \cite{K2}.

\end{ex}

\begin{ex}
(del Pezzo surfaces of degree 1)
This case is a degenerate case of the previous example.  Let $S$ be a smooth del Pezzo surface of degree 1.
The anti-bi-canonical model of $S$ is a double cover of a quadric cone $Q_0$ in $\bbP^3$ branched along 
the vertex of $Q_0$ and a smooth curve $C$ of genus 4.  It is known that $C$ is an intersection of $Q_0$ and a cubic surface.  By taking the triple cover of
$Q_0$ branched along $C$ and then taking its minimal resolution, we have a K3 surface $X$ with an automorphism of
order 3.  In this case, $M = U\oplus A_2(2)$ and $N = U \oplus U \oplus E_8 \oplus D_4\oplus A_2$.
Thus we have an isomorphism between the moduli space of smooth del Pezzo surfaces of degree 1 and 
$\tilde{\Gamma}_{M,\rho}\backslash (\bbB_8\setminus \Delta)$. For more details we refer the reader to the paper \cite{K2}.

\end{ex}

We summarize these examples in the following table.

\begin{table}[h]
\footnotesize{\[
\begin{array}{rllll}
{}& {\rm Deligne-Mostow}&M&N \\
{\rm Curves\ of\ genus\ 4} & ({1\over 6}{1\over 6}{1\over 6}{1\over 6}{1\over 6}{1\over 6}{1\over 6}{1\over 6}{1\over 6}{1\over 6}{1\over 6}{1\over 6})&U(3)& U\oplus U(3)  \oplus E_8\oplus E_8 \\
\noalign{\smallskip}
{\rm del\ Pezzo\ of\ degree\ 1} &{\rm a\ subball\ quotient} & U\oplus A_2(2)&U\oplus U  \oplus E_8\oplus D_4\oplus A_2  \\
    \noalign{\smallskip}
{\rm del\ Pezzo\ of\ degree\ 2} &{\rm not\ appear }& U(2) \oplus A_1^{\oplus 6} & U(2)\oplus U(2)  \oplus D_8 \oplus A_1^{\oplus 2}\\
    \noalign{\smallskip}
{\rm del\ Pezzo\ of\ degree\ 3} &({2\over 6}{2\over 6}{2\over 6}{2\over 6}{2\over 6}{1\over 6}{1\over 6}) & U \oplus A_2^{\oplus 5} & A_2(-1) \oplus  A_2^{\oplus 4} \\ 
\noalign{\smallskip}
{\rm del\ Pezzo\ of\ degree\ 4} & ({2\over 5}{2\over 5}{2\over 5}{2\over 5}{2\over 5}) & U \oplus D_8 \oplus D_8 & U(2) \oplus U(2)  \\ 
\noalign{\smallskip}
{\rm 6\ points\ on\ \bbP^1} & ({1\over 3}{1\over 3}{1\over 3}{1\over 3}{1\over 3}{1\over 3})& U \oplus E_6 \oplus A_2^{\oplus 3} & A_2(-1) \oplus A_2^{\oplus 3}\\
\noalign{\smallskip}
{\rm 8\ points\ on\ \bbP^1} & ({1\over 4}{1\over 4}{1\over 4}{1\over 4}{1\over 4}{1\over 4}{1\over 4}{1\over 4})& U(2) \oplus D_4 \oplus D_4 & U\oplus U(2) \oplus D_4 \oplus D_4\\
\end{array}
\]}
\end{table}

\section{Half-twists of Hodge structures} 

Here we discuss a version of a construction due to Bert van Geemen \cite{vG}. 

Let 
$$V_\bbC  = \bigoplus_{p+q=k} V^{pq}$$
be a polarized AHS of weight $k$ on a real vector space $V$. Let $\rho:A \to \GL(V)$ be an action of a cyclic group of order $d$ on $V$ by  Hodge isometries. Choose a generator $g$ of $A$ and  a subset $\Sigma$ of $\wt(\rho)$ which does not contain real characters and satisfies 
$\text{Im}(\chi(g)) > 0$ for any $\chi\in \Sigma$.
The vector space 
$$\bigoplus_{\chi\in \Sigma\cup \overline{\Sigma}}V_\bbC(\chi) \cong W_\bbC$$
for some  vector subspace $W\subset V$ (because it is invariant with respect to the conjugation of $V_\bbC$). 
Write 
\begin{equation}\label{twist1}
V_\Sigma^{pq} = \bigoplus_{\chi\in \Sigma} V_\chi^{pq},\quad V_{\overline{\Sigma}}^{pq} = \bigoplus_{\chi\in \overline{\Sigma}} V_\chi^{pq}.\end{equation}
Define the {\it negative half twist} of $V$ to be the decomposition 
$$W_\bbC = \bigoplus_{r+s = k+1} W^{rs},$$
where
$$W^{rs} = V_{\Sigma}^{r-1s}\oplus V_{\overline{\Sigma}}^{rs-1}.$$ 
Obviously, the decomposition \eqref{twist1} satisfies property (HD1) of AHS. Let us define a polarization form on $W$ by changing the polarization form $Q$ on $V$ with 
$$Q'(x,y) = Q(x,g(y))-Q(x,g^{-1}(y)).$$
Using that each $V_\bbC(\chi)$ and $V_\bbC(\chi')$ are orthogonal with respect to $Q$ unless $\bar{\chi} = \chi'$,
and $V^{ab}$ orthogonal to $V^{a'b'}$ unless $a=b'$, we check property (HD2).  
Let $x\in V_{\chi}^{r-1s}, \chi\in \Sigma$, write $\chi(g)-\bar{\chi}(g) = bi,$ where $b > 0$. We have
$$i^{r-s}Q'(x,\bar{x}) = i^{r-s}Q(x,g(\bar{x})-g^{-1}(\bar{x})) = $$
$$i^{r-s}Q(x,(\bar{\chi}(g)-\chi(g))\bar{x}) = i^{r-s}(\bar{\chi}(g)-\chi(g))Q(x,\bar{x}) = bi^{r-s-1}Q(x,\bar{x}) > 0.$$
Similarly, we check that $i^{r-s}Q'(x,\bar{x}) > 0$ if $x\in V_{\chi}^{rs-1}$, where $\chi\in \overline{\Sigma}$. This checks property (HD3). 

The situation with integral structure is more complicated. Let $F$ be the extension of $\bbQ$ obtained by joining a $d$th root of unity. Suppose our AHS has an integral structure with respect to some lattice $\Lambda$ in $V$ and $\rho$ is obtained from a representation in $\GL(\Lambda)$.  Then $\Lambda$ acquires a structure of a module over the ring of integers $\calO$ in the field  $F = \bbQ(e^{2\pi i/d})$ by setting
$e^{2\pi i/d}\cdot v = \rho(g)(v)$. The vector space $W$ becomes a vector space over $F$ and $\Lambda$ is a lattice in $V$ such that $Q'$ can be obtained from a bilinear form on $\Lambda$ taking its values in $\calO$. 

\begin{ex} Assume $V_\bbC = V^{10}\oplus V^{01}$ is a polarized AHS of weight 1. Then
 $$W^{20} = V_{\Sigma}^{10}, \quad W^{11} = V_{\Sigma}^{01}\oplus V_{\overline{\Sigma}}^{10}, \quad W^{02}  = V_{\overline{\Sigma}}^{01}.$$
Suppose $\dim V_{\Sigma}^{10} = 1$, then we obtain a Hodge structure of the same type as the Hodge structure arising from a K3 surface. It does not need to be a Hodge structure of a K3 surfaces if, say 
$\dim W^{11} > 19$. 

Let $\bldmu = (\mu_1,\ldots,\mu_m) = (\frac{a_1}{d},\ldots,\frac{a_m}{d})$ as in section \ref{sec7}. We assume that the condition \eqref{2} holds and $m > 4$. Consider the curve $X$ isomorphic to a nonsingular projective model of the affine  curve \eqref{eq}. Let $C_d$ be a cyclic group of order $d$ and $\chi\in \hat{A}$ be the character corresponding to $\bldmu$. We know that 
$\dim H^{10}(X,\bbC)(\chi)  = 1$ and 
$\dim H^{10}(X,\bbC)(\bar{\chi})  = m-3.$
Take $\Sigma = \{\chi,\bar{\chi}\}$.  Then we obtain that the half-twist $W$ is AHS with Hodge numbers $(1,2m-6,1)$ and admitting a Hodge isometry $g$. Assume that $W$  contains a lattice $N$ such that the AHS has an integral structure with respect to $N$. It is easy to see that the signature of $N$ is equal to $(2,2m-6)$. Also assume that $N$ can be primitively embedded in $L_{K3}$ and let $M = N^\perp$. Then our  AHS corresponds to a point in
$\calD_M^{\rho,\chi} \cong \bbB_{m-3}$. By surjectivity of the period map it corresponds to a $(M,\rho)$-polarized K3 surface. 
For $d=3,4,5,6$, these assumptions holds (see Section \ref{example}).  
It would be very interesting to study the remaining cases $d > 6$ and the monodromy group $\Gamma_\mu$ is an arithmetic subgroup of $\Aut(\bbB_{m-3}).$

\end{ex}


\end{document}